\documentclass[10pt]{article}
\usepackage{amsfonts, epsfig, amsmath, amssymb, color}

\usepackage[english]{babel}

\newcommand{\I}[1]{\mathbb{I}\{#1\} }
\newcommand{\E}{\mathbf{E}}
\renewcommand{\P}{\mathbf{P}}
\newcommand{\Px}{\mathbf{P}_{\!x}}

\renewcommand {\epsilon}{\varepsilon}

\newtheorem{theorem}{Theorem}[section]

\newtheorem{prop}{Proposition}[section]

\newtheorem{lemma}{Lemma}[section]

\newtheorem{rem}{Remark}[section]

\DeclareMathSymbol{\ophi}{\mathalpha}{letters}{"1E}

\newcommand{\e}{\varepsilon}
\renewcommand{\phi}{\varphi}

\newcommand{\be}{\begin{equation}}
\newcommand{\ee}{\end{equation}}
\newcommand{\ben}{\begin{equation*}}
\newcommand{\een}{\end{equation*}}

\newcommand{\ba}{\begin{equation}\begin{aligned}}
\newcommand{\ea}{\end{aligned}\end{equation}}

\newenvironment{proof}{\par\noindent{\bf Proof:}}{\hfill$\blacksquare$\par}
\newenvironment{Aproof}[2]{\par\noindent{\bf Proof of #1 #2:}}{\hfill$\blacksquare$\par}

\allowdisplaybreaks[4]

\begin{document}

\renewcommand{\figurename}{{\rm Fig.}}

\title{First exit times of solutions\\ of non-linear stochastic differential equations\\
driven by symmetric L\'evy processes\\ with $\alpha$-stable components\footnote{Last 
revision on 15 November 2005. 
To appear in \textit{Stochastic Processes and their Applications}
under the title \textit{First exit times of SDEs driven by stable L\'evy processes} (without figures).
}}
\date{September 2004}
\author{Peter Imkeller and Ilya Pavlyukevich\footnote{Corresponding author.}}

\maketitle

\abstract{
We study the exit problem of solutions of the stochastic differential equation
$dX_t^\e=-U'(X_t^\e)\,dt+\e\,dL_t$ from bounded or unbounded intervals which contain
the unique asymptotically stable critical point of the deterministic dynamical system
$\dot Y_t=-U'(Y_t)$. The process $L$ is composed of a standard Brownian motion and a
symmetric $\alpha$-stable L\'evy process. Using probabilistic estimates we
show that in the small noise limit $\e \to 0$, the exit time of $X^\e$ from an interval
is an exponentially distributed random variable and determine its expected value.
Due to the heavy-tail nature of the $\alpha$-stable component of $L$, the results
differ strongly from the well known case in which the deterministic dynamical system undergoes purely
Gaussian perturbations.
}



\section*{Introduction}

The study of dynamical systems subject to small random
perturbations keeps receiving much attention both in the physical and the
mathematical literature. In the simplest
(one-dimensional) setting, systems of this type find the following mathematical formulation.
Consider the ordinary
differential equation $\dot Y_t=-U'(Y_t)$, $Y_0\in [-b,a]$,
$a,b>0$, where $U$ is a potential function. Assume that $U(0)=0$ and
that $0$ is the unique asymptotically stable point of the deterministic dynamical
system associated with the equation, that means that for any starting point $x$ in $[-b, a]$,
the deterministic trajectory tends to 0: $Y_t(x)\to 0$ as
$t\to\infty$.\par\smallskip

Now perturb the deterministic dynamical system with some small random noise, that is,
consider the solutions of the stochastic differential equation
$$
X^\e_t=x-\int_0^t U'(X^\e_s)\,ds +\e \eta_t, \eqno (*)
$$
where $\eta$ is a random process, and the noise intensity
parameter $\e$ is small compared to the other parameters of the system ($\e$ is close to 0).
Under certain conditions on $U$ and $\eta$, for example under the assumptions that
$U'$ is Lipschitz and $\eta$ is a
semimartingale, the solution of the equation $(\ast)$ is well defined. In case $\eta$ is a standard
Brownian motion, the dynamical system is said to be perturbed by white noise. If $\eta$ is
an Ornstein-Uhlenbeck process, the terminology of `red noise' perturbation has been used.
The literature also knows perturbations by the so-called
shot noises, fractional Gaussian noises, or L\'evy noises.\par\smallskip

The stochastic dynamics of systems perturbed by white noise, which belong to the large class of
diffusions driven by the
Brownian motion, in the small noise limit, i.e.\ for $\e\to 0$, has received a great deal of attention
for decades and is particularly well understood. The pioneering papers on this topic are
\cite{PontryaginAV-33,Eyring-35,Kramers-40}.
Later on it was studied in
\cite{Friedman-74,Schuss-80,Williams-82,Day-83,BuslovM-92,KolokoltsovM-96,FreidlinW-98,Kolokoltsov-00,BovierEGK-04,BovierGK-05}, 
as well as in many other papers.\par\smallskip

One of the main results in this field is concerned with the time it takes for the diffusion to exit
a neighbourhood of a local attractor. Due to the fact that Kramers' pioneering paper
was one of the first to derive heuristically some properties of an \emph{exit law} of this type,
in particular in the physical literature it is often named Kramers' law.
Stated in modern language, it says that
the expected exit time is
exponentially large in $\e^{-2}$ and the growth rate can be interpreted as
the height of the potential
barrier to be overcome to leave the local attractor neighbourhood (see section \ref{s:heur} for
a rigorous formulation).\par\smallskip

White noise perturbations, however, are not always appropriate to
interpret real data in a reasonable way. This is the case for
example if the nature of the underlying random perturbation process
has to model abrupt pulses or \emph{extreme events}. A more natural
mathematical framework for these phenomena takes into account other
than purely Brownian perturbations. In particular infinitely
divisible L\'evy perturbations with jumps enter the
stage.\par\smallskip

The physical papers by P.~Ditlevsen
\cite{Ditlevsen-99b,Ditlevsen-99a} motivating our research
stipulate more general noise sources of the type
alluded to. They originate in simple physical concepts serving to interpret paleoclimatic data.
In fact, paleoclimatic records from the Greenland ice-core show that the climate of the
last glacial period experienced rapid transitions between cold basic glacial periods
and several warmer interstadials (the so-called Dansgaard--Oeschger events). Those records
are given by the concentration of certain oxygen, hydrogen or calcium isotopes in the
annual layers of the ice-core extending over several hundred millennia in the past.
They can be used to reconstruct the global Earth temperature for the time span for which the records
are available. The calcium signal has the highest --- almost annual ---
temporal resolution and provides the most conclusive information about the statistics
of the Dansgaard--Oeschger warming events. They start with a very rapid warming of the
North Atlantic region of about 5--10$^\circ$C within at most a few decades. The warming is
followed by a plateau phase with slow cooling extending over several centuries, ending with
an equally abrupt drop to basic glacial conditions.
The nature of these events is not clear, and several conceptual explanations have been
proposed. One line of arguments invokes
instabilities in the North Atlantic thermohaline circulation as a causal mechanism for these
abrupt climate changes, and for the millennial time scale between jumps. A different reasoning
starts from the idea that the coupled atmosphere-ocean system in the tropics possesses several meta-stable states,
and claims global teleconnections that trigger changes in the North Atlantic thermohaline circulation.
In \cite{GanopolskiR-01}, the effect of stochastic resonance was
brought into play in order to explain the observed random periodicity of the
Dansgaard--Oeschger events.\par\smallskip

As for many climate phenomena, due to the non-linearity of the system
the physical background is highly complex. To understand some basic features,
simple low-dimensional models may be used. For the phenomena under consideration,
in this spirit one may hope to recover important aspects of the statistical properties of the
observed data by modelling the paleoclimatic temperature process as the solution
of an equation of the type $(\ast)$ with some noise term whose nature has
to be determined. To take account of meta- or multi-stability, it is natural
to assume that
the climatic potential function has (at least) two wells, their stable minima corresponding
to one cold (basic glacial)
and at least one warm state (plus an intermediate one).
In this setting characteristics of the transition mechanism between climate states may
be reformulated in terms of the
exit problem from local attractor neighbourhoods for solutions of the stochastic differential equation.
This approach was taken in \cite{Ditlevsen-99b,Ditlevsen-99a}.
To account for a reasonable choice of random noise perturbing the system, a spectral analysis
of real ice-core data was performed in
\cite{Ditlevsen-99a}.
The obtained spectral decomposition exhibits a strong $\alpha$-stable
component with $\alpha\approx 1.75$. The paper \cite{Ditlevsen-99b}
is concerned with an analysis of the
exit times of $(\ast)$ with an $\alpha$-stable noise $\eta$ and in the limit of small
$\e$, performed on a physical level of rigour, with the help of a
fractional Fokker-Planck equation.\par\smallskip

Climate dynamics is not the only source of stochastic models in which
$\alpha$-stable noise appears.
For example, it was shown in \cite{SokolovMB-97,BrockmannS-02},
that the thermally activated motion of the
test particle along a polymer in three-dimensional space is subject to $\frac{1}{2}$-stable
motion due to the polymer's self-intersections. In recent years, L\'evy noise
sources have been playing
an increasingly important role in models of financial markets
(see for example \cite{EberleinR-99}).\par\smallskip

In this paper we consider the equation $(\ast)$ driven by the
L\'evy process $L$ which is the sum of a standard Brownian motion
and an $\alpha$-stable L\'evy motion. Our approach of the
asymptotic laws in the small noise limit of exit times from
bounded intervals or intervals which are unbounded from one side
is purely probabilistic and completely avoids fractional
Fokker-Planck equations. We understand it as a first step towards
a complete understanding of transition patterns of L\'evy-driven
dynamical systems in bi- or multi-stable potentials. The
mathematical challenge consists in a large deviation analysis for
exit times replacing the classical theory of Freidlin-Wentzell
\cite{FreidlinW-98} for diffusions
with Brownian noise. We base it on a noise intensity dependent
decomposition of $L$ into a sum of two independent processes: a
compound Poisson with large jumps on the one hand, and a sum of
the Brownian motion and a L\'evy motion with small jumps on the
other hand. Given such a decomposition for small noise intensity
$\e$, the main idea of our analysis is to prove that
asymptotically exits from the considered domains are due to large
jumps of the first component, while the second component is not
able to perturb the deterministic trajectory of solutions of
$(\ast)$ without noise essentially. For this reason, the usual
picture of a particle that has to climb a potential well being
pushed by a Brownian motion in order to exit a domain, which
captures the system's behaviour for Gaussian noise, changes
drastically here. Instead of the height of the potential well, a
large jump to exit just takes note of the distance from the
domain's boundary. Pure horizontal distances replace geometric
quantities related to the potential in the large deviations'
estimates for exit times in the $\alpha$-stable L\'evy case. Also,
the mean values of the exit times change essentially in comparison
to the Gaussian noise setting: instead of Kramers' times we obtain
exit times of the order of $\e^{-\alpha}$ in the small noise
limit, i.e.\ times of \emph{polynomial} instead of
\emph{exponential} dependence on $\e$. Our approach can be
extended to more general heavy-tailed  L\'evy processes.
\par\smallskip

The structure of the paper is as follows.
In section \ref{s:heur} we give a heuristic discussion of the asymptotic exit law
based on the decomposition of the $\alpha$-stable L\'evy noise
perturbing our system, and the heuristic
picture that the small jump component does not essentially affect the
asymptotic behaviour. In section \ref{s:dev}
we underpin this heuristic picture with mathematical rigour in proving
that trajectories of the deterministic
system and the one in which only the small noise component
is admitted to perturb are asymptotically very close.
This crucial observation is used in section \ref{s:law:sigma} to
derive in a rather technical way upper and
lower estimates for the law of the exit time of a bounded interval.
In sections \ref{s:devU} and \ref{s:law:tau}
these results are transferred to the setting of intervals which are
bounded on one side. This requires the
possibility for the deterministic trajectory to return from $-\infty$ in finite time.

\medskip

\textbf{Acknowledgements.} This work was partially supported by the
DFG Research Center \textit{Mathematics for Key Technologies
(Matheon)} and the DFG research project \textit{Stochastic Dynamics
of Climate States}.

\section*{Preliminaries and notation}

On a filtered probability space
$(\Omega,\mathcal{F},(\mathcal{F}_t)_{t\geq 0}, \P)$ we consider a
stochastic differential equation driven by a L\'evy noise of
intensity $\e$: \be \label{eq:X} X_t^{\e}=x-\int_0^t
U'(X_{s-}^{\e})\,ds+\e L_t,\quad \e>0. \ee In general, a L\'evy
process is known to be a random process with independent and stationary increments,
which is continuous in probability and possesses rcll
paths. It is completely
determined by its one-dimensional distributions which are infinitely
divisible and characterised by the L\'evy-Hin\v{c}in formula. In
this paper we assume that
\be
\mathbf E e^{i\lambda L_1}=
\exp\left\lbrace -d\frac{\lambda^2}{2}+
\int_{\mathbb{R}\backslash\{0\}} (e^{i\lambda y}-1-i\lambda y \I{|y|<1})
\frac{dy}{|y|^{1+\alpha}}\right\rbrace,
\ee
that is $L$ is a sum of
a standard Brownian motion with variance $d\geq 0$ and
an independent $\alpha$-stable L\'evy motion with $0<\alpha<2$. More information on
L\'evy processes can be obtained from \cite{Bertoin-98,Sato-99}. Since a L\'evy process is a semimartingale, the
standard theory of stochastic integration applies to equation
\eqref{eq:X}, see
\cite{SamorodnitskyG-03,Protter-04} for more
details. Throughout this paper we assume that the underlying filtration
fulfils the \emph{usual
conditions} in the sense of \cite{Protter-04}, i.e.\ the
filtration $(\mathcal{F}_t)_{t\geq 0}$
consists of $\sigma$-algebras which are complete with respect to $\P$ and is
right-continuous.

The L\'evy measure of $L$ is given by $\nu(dy)=\frac{dy}{|y|^{1+\alpha}}$,
$y\neq 0$. It is heavy-tailed and has infinite mass for all $\alpha\in (0,2)$,
due to a strong intensity of small jumps.

We impose some geometric conditions on the potential function $U$. First, we assume that
$U$ has a `parabolic' shape with its non-degenerate global minimum at the origin, i.e.\
$U'(x)x\geq 0$, $U(0) = 0,$ $U'(x)=0$ iff $x=0$, and $U''(0)=M>0$. Further, to guarantee
the existence of a strong unique solution of \eqref{eq:X} on $\mathbb{R}$ we demand that
$U'$ is at least locally Lipschitz and increases faster than a linear function at $\pm \infty$
(see also \cite{SamorodnitskyG-03,Protter-04}).
Moreover, in order to obtain some fine small-noise approximations of $X^\e$ in
section~\ref{s:dev}, we need that
$U\in\mathcal{C}^3$ in some sufficiently large interval containing the origin.

We shall study the first exit problem for the process $X^\e$ from
bounded and unbounded
intervals in the small noise limit $\e\to 0.$ In fact, we consider two cases.

\textbf{(B)} Let $I=[-b,a]$, $a,b>0$, and define the first exit time from $I$ as
\be
\sigma(\e)=\inf \{t\geq 0\,:\,X_t^{\e}\notin [-b,a]\}.
\ee

\textbf{(U)} Let $J=(-\infty,a]$, $a>0,$ and assume that for some
$c_1,c_2>0$ the regularity condition
$U(x)=c_1|x|^{2+c_2}$, $x\to -\infty$,  holds (see Remark~\ref{rm:2} on page \pageref{rm:2}).
In this case we study the one-sided counterpart of $\sigma(\e)$
defined by
\be
\tau(\e)=\inf \{t\geq 0\,:\,X_t^{\e}> a\}.
\ee

We shall investigate the laws of $\sigma(\e)$ and $\tau(\e)$, in particular their mean
values, as $\e\to 0$.

As a main tool of our analysis, we decompose the L\'evy process $L$
into $\e$-dependent small and large jump components.
In mathematical terms, we represent the process $L$ at any time $t$ as
a sum of two independent processes
$L_t=\xi_t^\e+\eta_t^\e$, with characteristic functions
\ba
\mathbf E e^{i\lambda \xi^\e_1}&=
\exp\left\lbrace -d\frac{\lambda^2}{2}+
\int_{\mathbb R\backslash\{0\}} (e^{i\lambda y}-1-i\lambda y \I{|y|<1})
\I{|y|\leq\frac{1}{\sqrt{\e}}}\frac{dy}{|y|^{1+\alpha}}\right\rbrace,\\
\mathbf E e^{i\lambda \eta^\e_1}&=
\exp\left\lbrace
\int_{\mathbb R\backslash\{0\}} (e^{i\lambda y}-1)
\I{|y|>\frac{1}{\sqrt{\e}}}\frac{dy}{|y|^{1+\alpha}}\right\rbrace.
\ea
The L\'evy measures corresponding of the processes $\xi^\e$ and $\eta^\e$ are
\begin{align}
\nu_\xi^\e(\,\cdot\,)=\nu(\,\cdot\,\cap \{0<|y|\leq\frac{1}{\sqrt{\e}}\}),\quad
\nu_\eta^\e(\,\cdot\,)=\nu(\,\cdot\,\cap \{|y|>\frac{1}{\sqrt{\e}}\}).
\end{align}
The process $\xi^\e$ has an infinite L\'evy measure with support
 $[-\frac{1}{\sqrt{\e}},\frac{1}{\sqrt{\e}}]\backslash\{0\}$,
and makes infinitely many jumps on any time interval of positive length.
The absolute value of its
jumps does not exceed $1/\!\sqrt{\e}$. It will be explained later in Remark
\ref{rm:1} on page \pageref{rm:1} why the threshold $1/\!\sqrt{\e}$ is chosen.

The L\'evy measure $\nu_\eta^\e(\cdot)$ of $\eta^\e$ is finite.
Denote
\be
\beta_\e=\nu_\eta^\e(\mathbb R)=
\int_{\mathbb{R}\backslash [-\frac{1}{\sqrt{\e}},\frac{1}{\sqrt{\e}}]}
\frac{dy}{|y|^{1+\alpha}}=\frac{2}{\alpha}\e^{\alpha/2}.
\ee

Then, $\eta^\e$ is a compound Poisson process with intensity
$\beta_\e$, and jumps distributed according to the law
$\beta_\e^{-1}\nu_\eta^\e(\cdot)$.

Denote $\tau_k$, $k\geq 0$, the arrival times of the jumps of $\eta^\e$ with $\tau_0=0$.
Let $T_k=\tau_k-\tau_{k-1}$ denote the inter-jump periods, and
$W_k=\eta^\e_{\tau_k}-\eta^\e_{\tau_k-}$
the jump heights of $\eta^\e$.
Then, the three processes $(T_k)_{k\ge 1},$ $(W_k)_{k\ge 1}$, and $\xi^\e$ are independent. Moreover,
\begin{equation}
\begin{aligned}
&\P(T_k\geq u)= \int_u^\infty \beta_\e e^{-\beta_\e s}\,ds=e^{-\beta_\e u},
\quad u\geq 0,\\
&\E T_k=\frac{1}{\beta_\e}=\frac{\alpha}{2}\e^{-\alpha/2} \to\infty \text{ as }\e\to 0,\\
&\P(W_k\in A)= \frac{1}{\beta_\e}\int_A
\I{|y|>\frac{1}{\sqrt{\e}}}\frac{1}{|y|^{1+\alpha}}\,dy, \quad
\text{for any Borel set}\quad A\subseteq \mathbb{R}.
\end{aligned}
\end{equation}
Due to the strong Markov property, for any stopping time $\tau$
the process $\xi^\e_{t+\tau}-\xi^\e_\tau$, $t\geq 0,$ is
also a L\'evy process with the same law as $\xi^\e$.

For $k\geq 1$ consider the processes
\begin{equation}
\begin{aligned}
\xi^k_t&=\xi^\e_{t+\tau_{k-1}}-\xi^\e_{\tau_{k-1}},\\
x^k_t(x)&=x-\int_0^t U'(x^k_{s-})\, ds +\e\xi^k_t, \quad t\in[0,T_k].
\end{aligned}
\end{equation}

In our notation, for $x\in\mathbb R$,
\begin{equation}
\begin{aligned}
X_t^{\e}&=x^1_t(x)+\e W_1\I{t=T_1}, \quad t\in[0,T_1],\\
X_{t+\tau_1}^{\e}&=x^2_t( x^1_{\tau_1}+\e W_1)+\e W_2\I{t=T_2}, \quad t\in[0,T_2],\\
&\cdots\\
X_{t+\tau_{k-1}}^{\e}&=x^k_t( x^{k-1}_{\tau_{k-1}}+\e W_{k-1})+\e W_k\I{t=T_k}, \quad t\in[0,T_k].
\end{aligned}
\end{equation}

Finally, we denote by $Y(x)$ the deterministic function solving the non-perturbed version of (\ref{eq:X})
\begin{equation}
\begin{aligned}\label{deterministic}
Y_t(x)=x-\int_0^t U'(Y_s(x))\,ds,\quad x\in\mathbb R.
\end{aligned}
\end{equation}
\
\section{Heuristic derivation of the main result and \\comparison with Gaussian case
\label{s:heur}}

In this section we shall provide the skeleton of a heuristic derivation of our main result on the
asymptotic law of the exit time from a bounded
interval. In the subsequent two sections a rigorous underpinning of these arguments will
be given.

On the interval $[0, T_k]$ and for $x\in\mathbb R$ let us consider $Y(x)$ and $x^k(x)$. These
processes satisfy the equations
\begin{equation}
\begin{aligned}
x^k_t(x)&=x-\int_0^t U'(x^k_{s-}(x))\,ds+\e\xi^k_t,\\
Y_t(x)&=x-\int_0^t U'(Y_s(x))\,ds, \quad t\in[0,T_k],
\end{aligned}
\end{equation}
while the law of $T_k$ is an exponential with parameter $\beta_\e.$

The process $\xi^k$ being the part of $L$ with the small jumps,
our analysis will be based on comparisons of
the trajectories of $x^k$ and $Y$. If they are close, e.g.\ if for some $\gamma>0$,
$\Px(\sup_{0\leq s\leq T_k} |x^k_s-Y_s(x)|> \e^\gamma)$ is small enough
as $\e\to 0$, we can apply the following reasoning.

For any $x\in I$, the deterministic solution $Y(x)$ converges exponentially fast
to the stable attractor $0$.
Define the relaxation time $R(x,\e)$ the process $Y$ needs to reach
an $\e^\gamma$-neighbourhood of $0$
from an arbitrary point $x\in I$. Then, as a separation of variables argument
in (\ref{deterministic}) implies, for some $\mu_1 >0$,
\be
R(x,\e)\leq \max\{\int_{-b}^{-\e^\gamma}\frac{dy}{-U'(y)}, \int_{\e^\gamma}^a\frac{dy}{U'(y)}\}
\leq \mu_1 |\ln\e|, \quad 0<\e\leq\e_0.
\ee

Since $Y(x)$ and $x^k(x)$ are close on the interjump interval
$[0,T_k]$ for all $k\geq 0$ and all $x\in I$, $X^\e$ can leave $I$
only at one of the time instants $\tau_k$ while jumping by the distance $\e W_k$.

\begin{figure}[h]
\begin{center}
\input{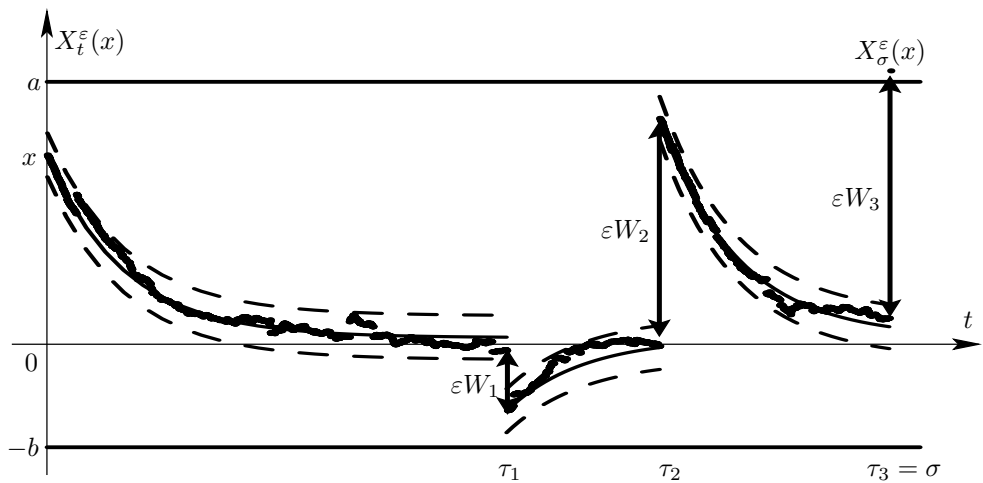}
\end{center}
\caption{A sample solution of the stochastic differential equation \eqref{eq:X}.\label{fig:2}}
\end{figure}

If the process has not left $I$ at jump number $k-1$, it waits for the next
possibility to jump at the end of a random exponentially distributed time period $T_k$.
Since $T_k$ (on average)
is essentially larger than the bound on the `relaxation' time
$\mu_1|\ln\e|$, $\e\to 0$, this means that $X^\e$ jumps from a small
neighbourhood of the attractor $0$ (see Fig.~\ref{fig:2}).

Therefore, the following heuristic estimate makes clear, that in the small noise limit
$\sigma(\e)$ is an exponentially
distributed random variable with parameter $\e^\alpha\theta/\alpha$, with
\be
\label{eq:theta}
\theta=\frac{1}{a^\alpha}+\frac{1}{b^\alpha}.
\ee
Indeed, for all
$k\geq 1$, $\tau_k=\sum_{j=1}^k T_j$ has a Gamma law with the density
$\beta_\e e^{-\beta_\e t}\frac{(\beta_\e t)^{k-1}}{(k-1)!}$. Hence for $u\geq 0$
\ba
\label{eq:fact}
\Px(\sigma(\e)>u)&\approx \sum_{k=1}^\infty \P(\tau_k>u)\cdot\Px(\sigma(\e)=\tau_k)\\
&=\sum_{k=1}^\infty \P(\tau_k>u) \cdot
\P(\e W_1\in I,\dots, \e W_{k-1}\in I,\e W_k\notin I)\\
&=\sum_{k=1}^\infty \int_u^\infty\beta_\e e^{-\beta_\e t}\frac{(\beta_\e t)^{k-1}}{(k-1)!}\, dt \cdot
(1-\P(\e W_1\notin I))^{k-1}\cdot\P(\e W_1\notin I)\\
&=\beta_\e\P(\e W_1\notin I) \int_u^\infty  e^{-\beta_\e t}
\sum_{k=1}^\infty \frac{ (\beta_\e t)^{k-1}(1-\P(\e W_1\notin I))^{k-1}}{(k-1)!}\, dt\\
&=\beta_\e\P(\e W_1\notin I)\int_u^\infty  e^{-\beta_\e t}  e^{\beta_\e t (1-\P(\e W_1\notin I))}\, dt
=e^{-u \beta_\e \P(\e W_1\notin I)}=e^{-u \e^\alpha\theta/\alpha}
\ea
For deriving the last formula we use the equations
\begin{equation}
\begin{aligned}
\P(\e W_1\notin I)&=\P(W_1< -\frac{b}{\e}\text{ or }W_1> \frac{a}{\e} )
=\frac{1}{\beta_\e}\left( \int_{\frac{a}{\e}}^\infty \frac{dy}{y^{1+\alpha}}
+\int_{\frac{b}{\e}}^\infty \frac{dy}{y^{1+\alpha}}\right) \\
&=\frac{1}{\beta_\e} \frac{1}{\alpha}
\left[ \left( \frac{\e}{a}\right)^\alpha+\left( \frac{\e}{b}\right)^\alpha\right] =
\frac{\e^\alpha}{\beta_\e\alpha}\theta.
\end{aligned}
\end{equation}
The mean value of $\sigma(\e)$ may be obtained immediately from \eqref{eq:fact},
or independently by the following reasoning:
\ba
\E_x\sigma(\e)&\approx \sum_{k=1}^{\infty} \E \tau_k\cdot \Px(\sigma(\e)=\tau_k)
=\sum_{k=1}^{\infty}k\E T_k\cdot\P(\e W_1\in I,\dots, \e W_{k-1}\in I, \e W_k\notin I)\\
&=\beta_\e\P(\e W_1\notin I)\sum_{k=1}^{\infty} k(1-\P(\e W_1\notin I))^{k-1}
=\frac{\beta_\e}{\P(\e W_1\notin I)}=
\frac{\alpha} {\e^\alpha}\left[ \frac{1}{a^\alpha}+\frac{1}{b^\alpha}\right]^{-1}.
\ea

The aim of this paper is to make these heuristic arguments rigorous.
This is done by proving the following Theorems.
\begin{theorem}
\label{th:main}
There exist positive constants $\e_0$, $\gamma$, $\delta$, and $C>0$ such that for
$0<\e\leq \e_0$
the following asymptotics holds
\begin{equation}
\exp\left\lbrace
-u \e^\alpha\frac{\theta}{\alpha} (1+C\e^\delta)
\right\rbrace
(1-C\e^\delta)\leq\Px(\sigma(\e)>u)\leq
\exp\left\lbrace
-u \e^\alpha\frac{\theta}{\alpha}(1-C\e^\delta)
\right\rbrace
(1+C\e^\delta)
\end{equation}
uniformly for all $x\in [-b+\e^\gamma, a-\e^\gamma]$ and $u\geq 0$, where $\theta = \frac{1}{a^{\alpha}} + \frac{1}{b^{\alpha}}$.
Consequently,
\begin{equation}
\label{eq:L1}
\E_x\sigma(\e)= \frac{\alpha}{\e^\alpha}\left[
\frac{1}{a^\alpha}+\frac{1}{b^\alpha}\right]^{-1} (1+\mathcal{O}(\e^\delta))
\end{equation}
uniformly for all $x\in [-b+\e^\gamma, a-\e^\gamma]$.
\end{theorem}

\begin{theorem}
\label{th:mainU}
There exist positive constants $\e_0$, $\gamma$, $\delta$, and $C>0$ such that for
$0<\e\leq \e_0$
the following asymptotics holds
\be
\exp\left\lbrace
-u \frac{\e^\alpha}{\alpha a^\alpha} (1+C\e^\delta)
\right\rbrace
(1-C\e^\delta)\leq\Px(\tau(\e)>u)\leq
\exp\left\lbrace
-u \frac{\e^\alpha}{\alpha a^\alpha}(1-C\e^\delta)
\right\rbrace
(1+C\e^\delta)
\ee
uniformly for all $x\leq a-\e^\gamma$ and $u\geq 0$.
Consequently,
\begin{equation}
\E_x\tau(\e)= \alpha \frac{a^\alpha}{\e^\alpha}
(1+\mathcal{O}(\e^\delta))
\end{equation}
uniformly for all $x\leq a-\e^\gamma$.
\end{theorem}

\vspace{.5cm}

It is interesting to compare the results stated above with their well-known counterparts
for diffusions driven by the Brownian motion of small intensity $\e$.
Together with \eqref{eq:X} consider the diffusion $\hat{X}^\e$ which solves the
stochastic differential equation
\be
\hat{X}^\e_t=x-\int_0^t U'(\hat{X}^\e_s)\, ds + \e W_t,
\ee
where $W$ is a standard one-dimensional Brownian motion, and $U$ is the same potential
as in \eqref{eq:X}.
For the diffusion $\hat{X}^\e$ we define the first exit time of the interval $I$ by
\be
\hat{\sigma}(\e)=\inf\{t\geq 0\,:\,\hat{X}^\e_t\notin [-b,a]\}.
\ee
Then the following statements hold for $\hat{\sigma}(\e)$ in the limit of small $\e$.

1.\ The first exit time $\hat{\sigma}(\e)$ is exponentially large in $\e^{-2}$. Assume for
definiteness, that $U(a)<U(-b)$. Then for any $\delta>0, x\in I,$ according to
\cite{FreidlinW-98}:
\be
%\label{eq:B1}
\Px(e^{(2U(a)-\delta)/\e^2}<\hat{\sigma}(\e)<e^{(2U(a)+\delta)/\e^2})\to 1\quad\text{as}
\quad \e\to 0.
\ee
Moreover, $\e^2\ln\E_x \hat{\sigma}(\e)\to 2U(a)$.

The mean of the first exit time can be calculated more explicitly (Kramers' law)
\cite{Kramers-40,Schuss-80,BovierEGK-04}:
\be
\label{eq:B2}
\E_x \hat{\sigma}(\e)\approx \frac{\e\sqrt{\pi}}{U'(a)\sqrt{U''(0)}}e^{2U(a)/\e^2}
\ee

2.\ The normalised first exit time is exponentially distributed
\cite{Williams-82, Day-83,BovierGK-05}: for $u\ge 0$
\be
%\label{eq:B3}
\Px\left(\frac{\hat{\sigma}(\e)}{\E_x \hat{\sigma}(\e)}>u\right)\to e^{-u}\quad\text{as}\quad \e\to 0,
\ee
uniformly in $x$ on compact subsets of $(-b,a)$.

As we see, $\hat{\sigma}(\e)$ and $\sigma(\e)$ have different orders of growth as $\e\to 0$.
The exit times of the $\alpha$-stable driven processes are much shorter because of the
presence of large jumps which occur with polynomially small probability.
To leave the interval, the diffusion $\hat{X}^\e$ has to overcome a potential barrier
of height either $U(-b)$ or $U(a)$. So in the case considered here,
$\hat{X}^\e_{\hat{\sigma}(\e)}=a$ with an overwhelming probability. The diffusion
`climbs' up in the potential landscape. This also explains why the pre-factor in \eqref{eq:B2}
depends on geometric properties of $U$ such as the slope at the exit point and the curvature
at the local minimum, the place where the diffusion spends most of its time before exit.

The process $X^\e$ on the contrary uses the possibility to exit the interval at one large jump.
This is the reason why the asymptotic exit time depends mainly on the distance between the stable point $0$
and the interval's boundaries.
The potential's geometry does not play a big role for the low order approximations of the
exit time $\sigma(\e)$. Although it is important for the proof,
it does not appear in the pre-factors of the mean exit time in \eqref{eq:L1} and remains
hidden in the error terms.

In the purely Gaussian case, to obtain the law of $\hat{\sigma}(\e)$, the theory of partial
differential equations is used. In fact, the probability $p_\e(x,u)=\Px(\hat{\sigma}(\e)>u)$ as
a function of $x$ and $u$
satisfies a backward Kolmogorov equation (parabolic partial differential equation)
with appropriate boundary conditions. The function $p_\e(x,u)$ can be
(at least in $\mathbb{R}$) expanded
in a Fourier series with respect to the eigensystem of the
diffusion's infinitesimal generator. Then, one concludes that
$p_\e(x,u)\approx e^{-\lambda_1^\e u}$, where $\lambda_1^\e$ is the first
eigenvalue. Further one shows, that
$\lambda_1^\e \E_x\hat{\sigma}(\e)\to 1$ as $\e\to 0$.
In contrast to this, in the present paper
we obtain results
without any use of operator theory.
This suggests that some asymptotic spectral properties of the integro-differential
operators corresponding to the process $X^\e$ can be formulated from the probabilistic
estimates obtained here. This is the subject of future research.

Perturbations of the deterministic dynamical systems by small infinitely divisible noises
were considered e.g.\ in \cite{FreidlinW-98},
however in a different setting. There, the small parameter $\e$ was  responsible for the
\emph{simultaneous} scaling of jump size and jump intensity. As the simplest example
of such a perturbation one can consider a compensated Poisson process with jump size of the
order $\e$ and jump intensity of the order $1/\e$ (see also
\cite{Borovkov-67}). In such a case the
dynamics in the limit corresponds to the one of the system perturbed by
white noise, i.e.\ the  probabilities  of the rare (exit) events are exponentially small in
$\e$ and the characteristic time scales are exponentially large.
Note that the perturbations considered in the present paper are of quite a different nature.
We only scale the jump sizes, while the jump intensities stay unchanged.

\section{Deviations from the deterministic trajectory: exit from bounded interval
\label{s:dev}}

In this section we estimate the deviation of the solutions of the stochastic differential
equation driven by the small-jump process $\e\xi^\e$ from the deterministic trajectory
on random time intervals of
exponentially distributed length. We show that the probabilities for at least polynomially small deviations
are polynomially small in $\e$ in the small noise limit. This rigorously underpins the starting point of our
heuristic derivation of the exit law.

For $x\in [-b,a]$ consider solutions $x^\e$ and $Y$ of the equations
\ba
\label{eq:eqs}
x_t^\e(x)&=x-\int_0^t U'(x_{s-}^\e(x))\, ds+\e\xi_t^\e,\\
Y_t(x)&=x-\int_0^t U'(Y_s(x))\, ds.
\ea
The goal of this section is to prove the following estimate.

\begin{prop}
\label{p:1}
Let $T_\e$ be exponentially distributed with parameter $\beta_\e$, and independent of
$\xi^\e$. Let $c>0$,
$\gamma=\frac{2-\alpha}{5}$.
Then there is $\e_0>0$ and $C>0$ such that for all
$0<\e\leq \e_0$ and $x\in [-b,a]$ the inequality
\be
\Px(\sup_{t\in [0,T_\e]}|x_t^\e-Y_t(x)|\geq c\e^\gamma)\leq C \e^{(\alpha+\gamma)/2}
\ee
holds.
\end{prop}

In order to prove the Proposition, we shall make use of the following Lemma in which
the estimation of the deviation
from the deterministic trajectory is executed on a bounded deterministic time interval.

\begin{lemma}
\label{l:1}
Let $T\geq 0$, $c>0$ and $\gamma=\frac{2-\alpha}{5}$ ($0<\gamma<\frac{2}{5}$). Then there exist
positive numbers $\e_0$ and $C$ such that for all
$0<\e\leq \e_0$ and $x\in [-b,a]$ the inequality
\be
\Px(\sup_{[0,T]}|x^\e_t-Y_t(x)|\geq c\e^\gamma)\leq C T \e^{\alpha+\gamma/2}
\ee
holds.
\end{lemma}

\begin{Aproof}{Proposition}{\ref{p:1}}
We apply Lemma \ref{l:1} and the definition of $\beta_\e$ to obtain for all $x\in [-b,a]$ and
$\e\leq \e_0$
\ba
\Px(\sup_{[0,T_\e]}|x^\e_t-Y_t(x)|\geq c\e^\gamma)&=
\int_0^\infty
\Px(\sup_{[0,\tau]}|x^\e_t-Y_t(x)|\geq c\e^\gamma)\beta_\e e^{-\beta_\e\tau}\, d\tau\\
&\leq C'\e^{\alpha+\gamma/2}\int_0^\infty \tau \beta_\e e^{-\beta_\e\tau}\, d\tau
\leq C\e^{(\alpha+\gamma)/2}.
\ea

\end{Aproof}

The proof of Lemma \ref{l:1} is performed in three Lemmas in the sequel.
It extensively uses the following geometric properties of
the potential $U$:
\begin{enumerate}
\item The deterministic trajectories $Y_t(x)$, $x\in [-b,a]$ converge to zero
as $t\to\infty$
due to the property $U'(x)x>0$ for $x\neq 0$.
\item The curvature of the potential at $x=0$ is positive. In small
neighbourhoods of zero we have
$U(x)= \frac{M}{2}x^2+o(x^2)$. Consequently $Y$ decays
there like
$e^{-Mt}$, and the dynamics of $x^\e$ reminds of the dynamics of a
process of Ornstein-Uhlenbeck type.
\end{enumerate}

We now prepare our rigorous analysis by an asymptotic expansion of the random
trajectories of $x^\e$ around the deterministic one of $Y$.
To this end, fix some $\delta>0$ small enough which will be specified later and define
\be
\hat{T}=\max\{\int_{-b}^{-\delta}\frac{dy}{-U'(y)}, \int_\delta^a\frac{dy}{U'(y)}\}.
\ee
$\hat{T}$ has the property that for all $x\in [-a,b]$ and $t\geq \hat{T}$,
$|Y_t(x)|\leq \delta$, i.e.\ after $\hat{T}$ the trajectory of $Y(x)$ is within a
$\delta$-neighbourhood of
the origin.
We next consider the representation of the process
$x^\e$ in powers of $\e$ of the form
\be
x^\e(x)=Y(x)+\e Z^\e(x)+R^\e(x),
\ee
where $Z^\e$ is the first approximation of $x^\e$ in powers of $\e$ satisfying the
stochastic differential equation
\be
\label{eq:Z}
Z_t^\e(x)=-\int_0^t U''(Y_s(x)) Z_{s-}^\e(x)\, ds+\xi_t^\e.
\ee
The solution to this equation is explicitly given by
\be
Z_t^\e(x)=\int_0^t e^{-\int_s^t U''(Y_u(x))\, du}\,d\xi_s^\e.
\ee
Integration by parts results in the following representation for $Z^\e$:
\be
Z_t^\e(x)=\xi_t^\e-
\int_0^t \xi_{s-}^\e U''(Y_s(x))e^{-\int_s^t U''(Y_u(x))\, du} \,ds.
\ee
For $x=0$, $Y_t(x)=0$ for all $t\geq 0$, and $Z^\e(0)$ is a process of the Ornstein-Uhlenbeck
type starting at zero and given by the equation
\be
Z_t^\e(0)=\xi_t^\e-
M\int_0^t \xi_{s-}^\e e^{-M(t-s)} \,ds.
\ee

\begin{lemma}
\label{l:Z} There is a universal constant $C_Z>2$ such that for
any $T>0$, $x\in [-b,a]$ and $\e>0$ 
\be
\sup_{s\in[0,T]}|Z^\e_s(x)|\leq C_Z \sup_{s\in[0,T]}|\xi^\e_s|, \quad \Px\text{-a.s.} 
\ee
\end{lemma}
\begin{proof}
Obviously, for $t\leq T$
\ba
|Z_t^\e(x)|&\leq\sup_{t\in [0,T]}|\xi_t^\e|
\left( 1+\sup_{t\in [0,T]} \int_0^t |U''(Y_s(x))| e^{-\int_s^t U''(Y_u(x))\, du} \,ds\right).
\ea
We shall show that the integral in the parenthesis is uniformly bounded.
Fix some $\delta>0$ small enough such that for
some $0<m_1<m_2$ the inequality
$m_1 < \inf_{|x|<\delta} U''(x)\leq \sup_{|x|<\delta} U''(x)< m_2$ holds. This implies, that
$m_1 < U''(Y_t(x))< m_2$ for all $x\in I$, and $t\geq \hat{T}$.
Let
\be
C_1=\max_{x\in I}
\int_0^{\hat{T}} |U''(Y_s(x))|e^{-\int_s^{\hat{T}} U''(Y_u(x))\, du} \,ds.
\ee
Consider an arbitrary $t\geq \hat{T}$. Then
\ba
\label{eq:2summand}
\int_0^t |U''(Y_s(x))&|e^{-\int_s^t U''(Y_u(x))\, du} \,ds\\
&=\int_0^{\hat{T}} |U''(Y_s(x))|e^{-\int_s^t U''(Y_u(x))\, du} \,ds
+\int_{\hat{T}}^t |U''(Y_s(x))|e^{-\int_s^t U''(Y_u(x))\, du} \,ds.
\ea
Let us estimate the first summand in \eqref{eq:2summand}. We have
\ba
\int_0^{\hat{T}} |U''(Y_s(x))|e^{-\int_s^t U''(Y_u(x))\, du} \,ds
&=e^{-\int_{\hat{T}}^t U''(Y_u(x))\, du}
\int_0^{\hat{T}} |U''(Y_s(x))|e^{-\int_s^{\hat{T}} U''(Y_u(x))\, du} \,ds\\
&\leq e^{-m_1(t-\hat{T})}C_1\leq C_1.
\ea
The second summand in \eqref{eq:2summand} is estimated analogously:
\ba
\int_{\hat{T}}^t |U''(Y_s(x))|e^{-\int_s^t U''(Y_u(x))\, du} \,ds\leq
m_2\int_{\hat{T}}^t e^{-m_1(t-s)}\,ds\leq \frac{m_2}{m_1}.
\ea
Taking $C_Z=\max\{2, C_1+\frac{m_2}{m_1}\}$ completes the proof.

\end{proof}

To estimate the remainder term $R^\e$ we need finer smoothness
properties of the potential $U$. However, the following Lemma shows that this restriction only has
to hold locally.

\begin{lemma}
\label{l:Rrough}
There exists $C>0$ such that for all $x\in[-b,a]$ and $T>0$,
\be
\sup_{t\in[0,T]}|R^\e_t(x)|\leq C
\ee
a.s.\ on the event $\{\omega \,:\, \sup_{t\in [0,T]}|\e\xi^\e_t(\omega)| <1\}$.

\end{lemma}
\begin{proof}
By hypothesis we know that for any $t\ge 0, x\in [-b,a]$ we have $-b\le Y_t(x) \le a.$
Moreover, on $\{\sup_{t\in [0,T]}|\e \xi_t^\e| < 1\}$ we have
$\sup_{t\in [0,T]}|\e Z_t^\e| < C_Z $ due to Lemma \ref{l:Z}.
Recall that $U'$ increases at least linearly at infinity. This guarantees the existence of
$C > 0$ such that for any $y\in [-b , a ], z\in [-C_Z, C_Z]$  we have
$$-U'(y+z+C)+ U'(y)+ U''(y)z<0.$$
Hence for any $T>0, T\ge\tau>0, x\in[-b,a]$ the inequality
\be
-U'(Y_{\tau}(x)+ \e Z_{\tau-}^\e(x)+C)+ U'(Y_{\tau}(x))+ U''(Y_{\tau}(x))\e Z_{\tau-}^\e(x)<0.
\ee
holds on the event $\{\sup_{t\in [0,T]}|\e \xi_t^\e| < 1\}.$ Now assume there is some $x\in[-b,a]$,
and some (smallest) $\tau\in[0,T]$
such that $R^\e_\tau(x)=C$. Observe that the rest term satisfies the integral equation
\be
R^\e_t(x)=\int_{0}^t f(R^\e_s(x), Y_s(x), \e Z^\e_{s-}(x))\,ds
\ee
with the smooth integrand
\ben
f(R, Y, \e Z)=-U'(Y+\e Z+R)+U'(Y)+U''(Y)(\e Z).
\een
This implicitly says that $R^\e$ is an absolutely continuous function of time. By definition of $\tau$,
we have
$$ 0 \le D^+R_t^{\e}(x)\vert_{t=\tau} = -U'(Y_{\tau}(x)+ \e Z_{\tau-}^\e(x)+C)+ U'(Y_{\tau}(x))+ U''(Y_{\tau}(x))\e Z_{\tau-}^\e(x)<0,
$$
a contradiction. A similar reasoning applies under the assumption
$R^\e_\tau(x)= - C.$ This completes the proof.
\end{proof}

Lemma \ref{l:Rrough} has a very convenient consequence. It states precisely that the solution process
$x^\e$, with initial state confined to $[-b,a],$ stays bounded by a deterministic constant on sets
of the form $\{\omega \,:\, \sup_{t\in [0,T]}|\e\xi^\e_t(\omega)| <1 \}.$ Therefore, in the small noise limit,
only local properties of $U$ are relevant to our analysis.

We next obtain a finer estimate of the remainder term $R^\e$
on the time interval $[0,\hat{T}]$.

\begin{lemma}
There exits $C_{\hat{T}}>0$ such that for $0<T\leq \hat{T}$
\be
\sup_{s\in [0,T]}|R^\e_s(x)|\leq  C_{\hat{T}} (\sup_{s\in [0,T]}|\e\xi^\e_s|)^2,
\quad \Px\text{-a.s.},
\ee
on the event $\{\omega \,:\, \sup_{[0,T]}|\e\xi^\e_t| <1 \}$
uniformly for $x\in [-b,a]$ and $\e>0$.
\end{lemma}

\begin{proof}
Using Lemma \ref{l:Rrough}, choose $K>0$ such that on the event $\{\omega \,:\, \sup_{t\in [0,T]}|\e\xi^\e_t| <1 \}$ the
processes $x^\e(x), \e Z^{\e}(x), R^\e(x)$ take their values in $[-K, K]$ as long as time runs in $[0,T].$
For $t\leq T\leq \hat{T}$, the rest term $R^\e$ satisfies the following integral equation:
\ba
R^\e_t(x)&=\int_0^t
\left[-U'(Y_s(x)+\e Z_{s-}^\e(x)+R^\e_s(x))+U'(Y_s(x))+U''(Y_s(x))\e Z_{s-}^\e(x)\right]\, ds\\
&=-\int_0^t \left[U'(Y_s(x)+\e Z_{s-}^\e(x)+R^\e_s(x))-U'(Y_s(x)+\e Z_{s-}^\e(x)) \right]\, ds\\
&-\int_0^t \left[U'(Y_s(x)+\e Z_{s-}^\e(x))-U'(Y_s(x))-U''(Y_s(x))\e Z_{s-}^\e(x)  \right]\, ds\\
&=-\int_0^t U''(\theta^1_{s-})R^\e_s(x)  \, ds
-\int_0^t \tfrac{1}{2}U^{(3)}(\theta^2_{s-})(\e Z_{s-}^\e(x))^2  \, ds
\ea
with appropriate
$\theta^1_s$, $\theta^2_s\in [-K,K]$.
Thus
\be
|R_t^\e(x)|\leq \int_0^t L|R^\e_s(x)|  \, ds +\tfrac{1}{2}TLC^2_Z(\sup_{[0,T]}|\e\xi^\e_t|)^2.
\ee
An application of Gronwall's lemma yields the final estimates
\ba
|R_t^\e(x)|&\leq \tfrac{1}{2}TLC^2_Z e^{TL} (\sup_{t\in [0,T]}|\e\xi^\e_t|)^2
\leq \tfrac{1}{2}\hat{T}LC^2_Z e^{\hat{T}L} (\sup_{t\in [0,T]}|\e\xi^\e_t|)^2=
C_{\hat{T}} (\sup_{t\in [0,T]}|\e\xi^\e_t|)^2,\\
\sup_{t\in [0,T]}|R_t^\e(x)|&\leq C_{\hat{T}} (\sup_{t\in [0,T]}|\e\xi^\e_t|)^2.
\ea

\end{proof}

Now we derive a uniform estimate of the rest term $R^\e$ on time intervals
longer than $\hat{T}$.

\begin{lemma}
\label{l:R}
There exist positive constants $C_R$ and $C_E\leq 1$ such that for $T\geq 0$
\be
\sup_{s\in [0,T]}|R^\e_s(x)|\leq C_R (\sup_{s\in [0,T]}|\e\xi^\e_s|)^2
\ee
on the event
\be
E_T=\{ \omega\,:\,  \sup_{t\in [0,T]}|\e\xi^\e_t|<C_E \}
\ee
uniformly for $x\in [-b,a]$.
\end{lemma}

\begin{proof}
Fix some positive $T\geq \hat{T}$ and let
$\omega \in \{ \omega\,:\,  \sup_{t\in [0,T]}|\e\xi^\e_t|<1 \}$.
Again, using Lemma \ref{l:Rrough}, choose $K>0$ such that on the event $\{\omega \,:\, \sup_{t\in [0,T]}|\e\xi^\e_t| <1 \}$ the
processes $x^\e(x), \e Z^{\e}(x), R^\e(x)$ take their values in $[-K, K]$ as long as time runs in $[0,T].$
The rest term $R^\e$ satisfies the integral equation
\be
R^\e_t(x)=R^\e_{\hat{T}}(x)+\int_{\hat{T}}^t f(R^\e_s(x), Y_s(x), \e Z^\e_{s-}(x))\,ds
\ee
with
\ben
f(R, Y, \e Z)=-U'(Y+\e Z+R)+U'(Y)+U''(Y)(\e Z).
\een
Moreover, $R^\e$ is an absolutely continuous function of time.
We write the Taylor expansion for the integrand $f$ with some $\theta\in\mathbb [-K,K]$:
\ba
f(R, Y, \e Z)&=-U'(Y+\e Z+R)+U'(Y)+U''(Y)(\e Z)\\
&=-U'(Y)-U''(Y)(R+\e Z)-\frac{U^{(3)}(\theta)}{2}(R+\e Z)^2+U'(Y)+ U''(Y)(\e Z) \\
&=-U''(Y)R-\frac{U^{(3)}(\theta)}{2}(R+\e Z)^2
\ea
Since $U\in {\cal C}^3$, $|U^{(3)}|$ is bounded, say by $L$, on $[-K,K].$ Using the
inequality $(R+\e Z)^2\leq 2(R^2+\e^2 Z^2)$ we
obtain that for all $R$, $Y$ and $Z$ such that $L |\e Z|<A$
\ba
f(R,Y,\e Z)&\leq -U''(Y)R +L R^2 + L (\e Z)^2 < -U''(Y)R +L R^2 + A=g^+(R,Y)\\
f(R,Y,\e Z)&\geq -U''(Y)R -L R^2 - L (\e Z)^2 > -U''(Y)R -L R^2 - A=g^-(R,Y).
\ea
With the help of Lemma~\ref{l:Z} this immediately implies for $\hat{T}<t\leq T$ that
\ba
R^\e_t(x)&=R^\e_{\hat{T}}(x)+\int_{\hat{T}}^t f(R^\e_s(x), Y_s(x), \e Z^\e_{s-}(x))\,ds
<R^\e_{\hat{T}}(x)+\int_{\hat{T}}^t g^+(R^\e_s(x),Y_s(x))\,ds,\\
R^\e_t(x)&=R^\e_{\hat{T}}(x)+\int_{\hat{T}}^t f(R^\e_s(x), Y_s(x), \e Z^\e_{s-}(x))\,ds
>R^\e_{\hat{T}}(x)+\int_{\hat{T}}^t g^-(R^\e_s(x),Y_s(x))\,ds.
\ea
with $A=D \e^2  (\sup_{t\in [0,T]}|\xi^\e_t|)^2$ and a constant $D>2LC_Z^2$
which will be specified later.

To estimate the rest term $R^\e$ on $[\hat{T},T]$ we apply the subsequent comparison Lemma \ref{l:c2}.
Consider an event
\ba
E_1=\{\omega\,:\, |R^\e_{\hat{T}}(x)|<\frac{A}{m_2}\}
\supseteq
\{\omega\,:\, C_{\hat{T}}(\sup_{t\in [0,\hat{T}]}|\e\xi^\e_t|)^2<
\frac{D}{m_2}(\sup_{t\in [0,T]}|\e\xi^\e_t|)^2\}
\ea
Thus, setting $D>\max\{2LC_Z^2, C_{\hat{T}}m_2\}$ we obtain that $\Px(E_1)=1$, $x\in[-b,a]$
and the
conditions of Lemma  \ref{l:c2} are fulfilled on the event
\ba
E_2&=\{\omega\,:\, m_1^2-4AL > 0 \}=\{ \omega\,:\,  \sup_{t\in [0,T]}|\e\xi^\e_t|<\frac{m_1}{2\sqrt{LD}} \}.
\ea
On $E_2$ the following inequality holds a.s.\ for $x\in[-b,a]$:
\be
\sup_{t\in [\hat{T},T]}|R_t^\e(x)|\leq \frac{2D}{m_2}  (\sup_{t\in [0,T]}|\e\xi^\e_t|)^2.
\ee
Thus denoting $C_R=\max\{C_{\hat{T}},\frac{m_1}{2\sqrt{LD}}\}$ and
$C_E=\min\{\frac{2D}{m_2},1\}$ we may finish the proof.

\end{proof}

\begin{lemma}[Comparison lemma]
\label{l:c2}
Let $T\geq 0$, $Y$ be a smooth function on $[0,T]$ and $Z$ a rcll
function on $[0,T]$.
Consider the integral equation
\be
\label{eq:Rlemma}
R_t=R_0+\int_0^t f(R_s,Y_s,Z_{s-})\,ds,\quad t\in [0,T],\\
\ee
with a smooth function $f$ satisfying
\ba
g^-(R,Y)&<f(R,Y,Z)<g^+(R,Y), \quad R, Y, Z\in\mathbb R,\\
g^\pm(R,Y)&=-U''(Y)R \pm L R^2\pm A, \quad R, Y \in \mathbb{R},  \quad \text{with}\quad L,A>0.
\ea
Moreover, let $0<m_1<U''(Y_t)<m_2$, $t\in[0,T]$, and  $m_1^2-4AL>0$.
Then for $0<t\leq T$ the following holds:
\begin{enumerate}
\item if $R_0< \frac{A}{m_2}$ then $R_t<\frac{2A}{m_2}$;
\item if $R_0> -\frac{A}{m_1}$ then $R_t>-\frac{2A}{m_1}$.
\end{enumerate}
This yields, that
if $|R_0|<\frac{A}{m_2}$ then $\sup_{t\in [0,T]}|R_t|<\frac{2A}{m_2}$.
\end{lemma}

\begin{proof}
\noindent To prove the first statement, together with \eqref{eq:Rlemma} consider the Riccati equation
\be
r^+_t=R_0+\int_0^t g^+(r^+_s,Y_s)\,ds,\quad t\in[0,T].
\ee
Under the conditions of the Lemma, it is enough to prove two statements:
\begin{description}
\item[a)] $R_t<r^+_t$ for $0<t\leq T$.
\item[b)] $r^+_t<\frac{2A}{m_1}$ for $t\geq 0$.
\end{description}
To show \textbf{a)} we note that at the starting point $t=0$,
\be
D^+R_t\Big|_{t=0}=\lim_{h\downarrow 0}\frac{R_h-R_0}{h}
=f(R_0,Y_0,Z_0)<g^+(R_0,Y_0)=\dot r^+_t\Big|_{t=0},
\ee
consequently it follows from the continuity of $R$ and $r^+$ that $r^+_t>R_t$ for
at least positive and small $t$.
Assume there exists $\tau=\inf\{t>0\,:\, R_\tau=r^+_\tau\}$ such that $\tau\in (0,T]$.
At the point $\tau$ the left derivative of $R$ is necessarily not less than the
derivative of $r^+$ which leads to the following contradiction:
\ba
D^-R_t\Big|_{t=\tau}&=\lim_{h\downarrow 0}\frac{R_\tau-R_{\tau-h}}{h}
=f(R_\tau,Y_\tau,Z_{\tau-})\geq \dot r^+_t\Big|_{t=\tau}=g^+(r^+_\tau,Y_\tau),\\
f(R_\tau,Y_\tau,Z_{\tau-})&=f(r^+_\tau,Y_\tau,Z_{\tau-})<g^+(r_\tau,Y_\tau).
\ea
To prove \textbf{b)}, we compare $r^+$ with the stationary solution of the autonomous Riccati
equation
\be
\label{eq:p}
p_t=p_0+\int_0^t (-m_2 p_s+Lp_s^2+A)\, ds,\quad t\geq 0.
\ee
Equation  \eqref{eq:p} has two positive stationary solutions at which the integrand vanishes:
\be
p_0=p^\pm=\frac{m_2}{2L}\left(1\pm\sqrt{1-\frac{4AL}{m_2^2}} \right) .
\ee
Repeating the comparison argument used for \textbf{a)}, we obtain that
if $R_0<p^-$, then  $r^+_t<p^-$, $t\in[0,T]$.
Applying the elementary inequalities
\be
\frac{x}{2}\leq 1-\sqrt{1-x}\leq x, \quad x\in[0,1],
\ee
to $p^-$ we obtain that $\frac{A}{m_2}\leq p^-\leq \frac{2A}{m_2}<\frac{2A}{m_1}$.
This guarantees that for $R_0<\frac{A}{m_2}$, the solution of \eqref{eq:Rlemma}
does not exceed $\frac{2A}{m_1}$.

The proof of the second statement is analogous.
\end{proof}

\medskip

\begin{Aproof}{Lemma}{\ref{l:1}}
Let $T\geq 0$ and $x\in [-b,a]$. Choose $E_T$ according to Lemma \ref{l:R}.
Then there exists $\e_0$ such that for
$\e\leq\e_0$ the following holds:
\ba
&\{\sup_{t\in [0,T]}|x^\e_t(x)-Y_t(x)|\geq c\e^\gamma\}=
\{\sup_{t\in [0,T]}|\e Z^\e_t(x)+R^\e_t(x)|\geq c\e^\gamma\}\\
&\subseteq
\{\sup_{t\in [0,T]}|\e Z^\e_t(x)|\geq \tfrac{c\e^\gamma}{2}\}\cup
\{\sup_{t\in [0,T]}|R^\e_t(x)|\geq \tfrac{c\e^\gamma}{2}\}\\
&\subseteq
\{\sup_{t\in [0,T]}|\e\xi^\e_t|\geq \tfrac{c\e^\gamma}{2C_Z}\}\cup
\left[ \{\sup_{t\in [0,T]}|R^\e_t(x)|\geq \tfrac{c\e^\gamma}{2}\}\cap E_T \right] \cup
\left[\{\sup_{t\in [0,T]}|R^\e_t(x)|\geq \tfrac{c\e^\gamma}{2}\} \cap E^c_T\right]\\
&\subseteq
\{\sup_{t\in [0,T]}|\e\xi^\e_t|\geq \tfrac{c\e^\gamma}{2C_Z}\}\cup
\{\sup_{t\in [0,T]}|\e\xi^\e_t|\geq \tfrac{c\e^{\gamma/2}}{2\sqrt{C_R}}\}\cup
\{\sup_{t\in [0,T]}|\e\xi^\e_t|\geq C_E\}\\
&=
\{\sup_{t\in [0,T]}|\e\xi^\e_t|\geq \tfrac{c\e^\gamma}{2C_Z}\}.
\ea
Consequently, with the help of Kolmogorov's inequality we obtain for small $\e$ and some
$C>0$
\ba
\Px&(\sup_{t\in [0,T]}|x^\e_t-Y_t(x)|\geq c\e^\gamma)\leq
\P(\sup_{t\in [0,T]}|\e \xi^\e_t|\geq \frac{c}{2C_Z}\e^\gamma)\leq
\frac{4C_Z^2}{c^2\e^{2\gamma}}\E(\e\xi^\e_T)^2\\
&=T\frac{4C_Z^2}{c^2}(\frac{2}{2-\alpha}\e^{1+\alpha/2-2\gamma}+d\e^{2-2\gamma})
\leq C T\e^{1+\alpha/2-2\gamma}.
\ea
This completes the proof.

\end{Aproof}

\section{The law of $\sigma(\e)$\label{s:law:sigma}}

For the purposes of this rather technical section we introduce the following
notation.
Denote $W_0=T_0=0$, $x^1(0)=x$, and write $\mathbb I A$ for the indicator function of a measurable set $A$.
Let $I^\pm_q$ denote the exterior and interior neighbourhoods of $I$ defined by
$I^\pm_q=[-b\mp q, a\pm q]$ for $0\leq q< \min\{a,b\}$.

\begin{figure}
\begin{center}
\input{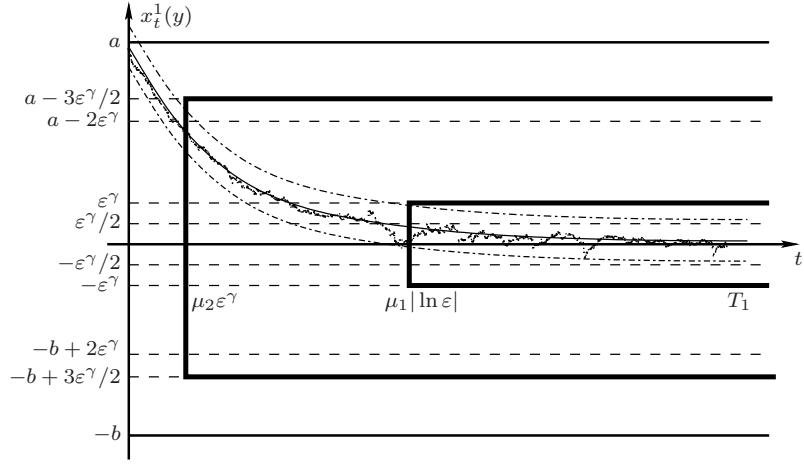}
\end{center}
\caption{The dynamics of $x^1(y)$ under condition
$\sup|x^1_t(y)-Y_t(y)|\leq \tfrac{\e^\gamma}{2}$.
\label{f:1} }
\end{figure}

Throughout this section we use the following constants.
Choose $\mu_1$, $\mu_2$ and $\e_0>0$ such that for all $y\in I$ the following holds for
$\e\in (0,\e_0]$, (see Fig.~\ref{f:1}):

\noindent
1.\ $|Y_t(y)|\leq \tfrac{\e^\gamma}{2}$ for  $t\geq \mu_1|\ln\e|$.
Consequently, denoting
$z_t^k(y) = x^k_t(y) - Y_t(y),$ we have $|x^k_t(y)|\leq\e^\gamma$ for $t\geq \mu_1|\ln\e|$
under the condition
$|z^k_t(y)| \leq \tfrac{\e^\gamma}{2}$.

\noindent
2.\ $Y_t(y)\in I^-_{2\e^\gamma}$ for all $t\geq \mu_2\e^\gamma$.
Consequently, $|x^k_t(y)|\in I^-_{3\e^\gamma/2}$ for $t\geq \mu_2\e^\gamma$ under the condition
$|z^k_t(y)| \leq \tfrac{\e^\gamma}{2}$.
Let us show that $\mu_2$ with this property exists. Indeed, for all $y\in I$, $|Y_t(y)|$ is
strictly decreasing
($y\neq 0$), and $Y_t(0)\equiv 0$. Moreover, for $-b\leq y_1<y_2\leq a$, $Y_t(y_1)<Y_t(y_2)$.
Let $r(x,\e)=\inf\{t\geq 0\,:\, Y_t(x)\in I^-_{2\e^\gamma}\}$. Then by comparison for all $x\in I$
\begin{equation}
r(x,\e)\leq \max\{\int^a_{a-2\e^\gamma} \frac{dy}{U'(y)},
\int^{-b+2\e^\gamma}_{-b} \frac{dy}{|U'(y)|}    \}\leq \mu_2\e^\gamma, \quad 0<\e\leq\e_0.
\end{equation}

%%%
%%%
\subsection{Upper estimate \label{S:Si:above}}
%%%
%%%

In this subsection we give estimates of $\Px(\sigma(\e)>u)$
from above as $\e\to 0$, $u>0$. They are comprised
in the following Proposition with a rather technical proof. Recall that $\gamma$
has been chosen according to
$\alpha$ in Proposition \ref{p:1}.

\begin{prop}
\label{p:above}
Let $\delta=\min\{\alpha/2,\gamma/2\}$. There exist constants $\e_0>0$ and $C>0$ such that for all $0<\e\leq \e_0,$
$x\in [-b+\e^\gamma, a-\e^\gamma]$ and $u\geq 0$ the following inequality holds
\be
\Px(\sigma(\e)>u)\leq
\exp\left\lbrace
-u \frac{\e^\alpha}{\alpha}\left[\frac{1}{a^\alpha}+\frac{1}{b^\alpha}\right] (1-C\e^\delta)
\right\rbrace
(1+C\e^\delta).
\ee
\end{prop}
\begin{proof}
For $x\in I$, we use the following obvious inequality
\ba
\Px(\sigma(\e)>u)&=\sum_{k=1}^\infty \P(\tau_k>u)\Px(\sigma(\e)=\tau_k)+
\Px(\sigma(\e)>u|\sigma(\e)\in (\tau_{k-1},\tau_k))\Px(\sigma(\e)\in (\tau_{k-1},\tau_k))\\
&\leq  \sum_{k=1}^\infty \P(\tau_k>u)\Big[ \Px(\sigma(\e)=\tau_k)+
\Px(\sigma(\e)\in (\tau_{k-1},\tau_k))\Big].
\ea

Then for any $x\in I$ and $k\in \mathbb N$, applying the independence and law properties of the
processes $x^j$, $j\in\mathbb N$, the following chain of inequalities is deduced which
results in a factorisation formula for the
probability under estimation (compare with \eqref{eq:fact}):
\begin{equation}
\label{eq:stauk}
\begin{aligned}
\Px(\sigma(\e)&=\tau_k)=\E_x\I{X_s^\e\in I, s\in[0,\tau_k), X_{\tau_k}^\e\notin I}\\
&= \E_x\prod_{j=1}^{k-1}
\I{x^j_s(x^{j-1}_{T_{j-1}}+\e W_{j-1})\in I, s\in[0,T_j],
x^j_{T_j}(x^{j-1}_{T_{j-1}}+\e W_{j-1})+\e W_j\in I} \\
&\qquad\qquad\times
\I{x^k_s(x^{k-1}_{T_{k-1}}+\e W_{k-1})\in I, s\in[0,T_k],
x^k_{T_k}(x^{k-1}_{T_{k-1}}+\e W_{k-1})+\e W_k\notin I}\\
&\leq \E\prod_{j=1}^{k-1}
\sup_{y\in I}\I{x^j_s(y)\in I, s\in[0,T_j], x^j_{T_j}(y)+\e W_j\in I} \\
&\qquad\qquad\times
\sup_{y\in I}\I{x^k_s(y)\in I, s\in[0,T_k], x^k_{T_k}(y)+\e W_k\notin I}\\
&=\prod_{j=1}^{k-1}
\E\left[ \sup_{y\in I}\I{x^j_s(y)\in I, s\in[0,T_j], x^j_{T_j}(y)+\e W_j\in I}\right] \\
&\qquad\qquad\times
\E \left[ \sup_{y\in I}\I{x^k_s(y)\in I, s\in[0,T_k],
x^k_{T_k}(y)+\e W_k\notin I}\right] \\
=&\left( \E\left[  \sup_{y\in I}\I{x^1_s(y)\in I, s\in[0,T_1],
x^1_{T_1}(y)+\e W_1\in I}\right]\right)^{k-1}  \\
&\qquad\qquad\times
\E \left[ \sup_{y\in I}\I{x^1_s(y)\in I, s\in[0,T_1],
x^1_{T_1}(y)+\e W_1\notin I}\right].
\end{aligned}
\end{equation}
Analogously we estimate the probability to exit between the $k$-th arrival times of
the compound Poisson process $\eta^\e$, $k\in\mathbb N$.
Here we distinguish two cases.
\ba
&\text{In the first case} \quad k=1,\-x\in I^-_{\e^\gamma}.\quad\text{Then}\\
&\Px(\sigma(\e)\in(\tau_0,\tau_1))=\Px(\sigma(\e)\in(0,T_1))=
\E_x\I{X_s^\e\notin I \text{ for some } s\in (0,T_1)}\\
&\leq
\E \left[ \sup_{y\in I^-_{\e^\gamma}}\I{x^1_s(y)\notin I\text{ for some } s\in[0,T_1]}\right],\\
%%%
&\text{In the second case} \quad k\geq 2, \-x\in I.\quad\text{Then}\\
&\Px(\sigma(\e)\in(\tau_{k-1},\tau_k))=
\E_x\I{X_s^\e\in I, s\in[0,\tau_{k-1}],
X_s^\e\notin I \text{ for some } s\in (\tau_{k-1},\tau_k)}\\
&= \E_x \prod_{j=1}^{k-1}
\I{x^j_s(x^{j-1}_{T_{j-1}}+\e W_{j-1})\in I, s\in[0,T_j],
x^j_{T_j}(x^{j-1}_{T_{j-1}}+\e W_{j-1})+\e W_j\in I} \\
&\qquad\qquad\qquad\qquad\times
\I{
x^k_s(x^{k-1}_{T_{k-1}}+\e W_{k-1})\notin I \text{ for some } s\in [0,T_k)}\\
&\leq \E\prod_{j=1}^{k-2}
\sup_{y\in I}\I{x^j_s(y)\in I, s\in(0,T_j], x^j_{T_j}(y)+\e W_j\in I}\times\\
&\times
\sup_{y\in I}
\I{x^{k-1}_s(y)\in I, s\in(0,T_{k-1}], x^{k-1}_{T_{k-1}}(y)+\e W_{k-1}\in I,
x^k_s(x^{k-1}_{T_{k-1}}(y)+\e W_{k-1})\notin I\text{ for some } s\in[0,T_k]}\\
&=\left( \E\left[  \sup_{y\in I}\I{x^1_s(y)\in I, s\in[0,T_1],
x^1_{T_1}(y)+\e W_1\in I}\right]\right) ^{k-2}  \\
&\times
\E \left[ \sup_{y\in I}
\I{x^1_s(y)\in I, s\in(0,T_1], x^1_{T_1}(y)+\e W_1\in I,
x^2_s(x^1_{T_1}(y)+\e W_1)\notin I\text{ for some } s\in[0,T_2]}\right].
\ea
Next we specify separately in four steps the further estimation for the four different
events appearing in the formulae for
$\Px(\sigma(\e)=\tau_k)$ and $\Px(\sigma(\e)\in(\tau_{k-1},\tau_k))$.

\noindent
\textbf{Step 1.} Consider  $\I{x^1_s(y)\in I, s\in [0,T_1], x^1_{T_1}(y)+\e W_1\in I}$.
For $y\in I$, we may estimate
\ba
\label{eq:a1}
&\I{x^1_s(y)\in I, s\in [0,T_1], x^1_{T_1}(y)+\e W_1\in I}\\
&=\I{x^1_s(y)\in I, s\in [0,T_1], x^1_{T_1}(y)+\e W_1\in I}
\left(\I{\sup_{s\in [0,T_1]}|z^1_s(y)|> \tfrac{\e^\gamma}{2}}
+\I{\sup_{s\in [0,T_1]}|z^1_s(y)|\leq \tfrac{\e^\gamma}{2}}\right) \\
&\leq
\I{\sup_{s\in [0,T_1]}|z^1_s(y)|> \tfrac{\e^\gamma}{2}}+
\I{\sup_{s\in [0,T_1]}|z^1_s(y)|\leq \tfrac{\e^\gamma}{2},
 x^1_{T_1}(y)+\e W_1\in I}\\
&=\I{\sup_{s\in [0,T_1]}|z^1_s(y)|> \tfrac{\e^\gamma}{2}}\\
&+\I{\sup_{s\in [0,T_1]}|z^1_s(y)|\leq \tfrac{\e^\gamma}{2},  x^1_{T_1}(y)+\e W_1\in I}
\Bigg( \I{|\e W_1|\leq\tfrac{\e^\gamma}{2}}+\I{|\e W_1|>\tfrac{\e^\gamma}{2}}\Bigg) \\
&\leq\I{\sup_{s\in [0,T_1]}|z^1_s(y)|> \tfrac{\e^\gamma}{2}}
+\I{|\e W_1|\leq\tfrac{\e^\gamma}{2}}\\
&+\I{\sup_{s\in [0,T_1]}|z^1_s(y)|\leq \tfrac{\e^\gamma}{2},
|\e W_1|>\tfrac{\e^\gamma}{2},  x^1_{T_1}(y)+ \e W_1 \in I, T_1\geq \mu_1|\ln\e|}\\
&+\I{\sup_{s\in [0,T_1]}|z^1_s(y)|\leq \tfrac{\e^\gamma}{2},
|\e W_1|>\tfrac{\e^\gamma}{2}, x^1_{T_1}(y)+\e W_1\in I, T_1< \mu_1|\ln\e|}\\
&\leq\I{\sup_{s\in [0,T_1]}|z^1_s(y)|> \tfrac{\e^\gamma}{2}}
+\I{|\e W_1|\leq\tfrac{\e^\gamma}{2}}\\
&+\I{|\e W_1|>\tfrac{\e^\gamma}{2},\e W_1 \in I^+_{\e^\gamma}}
+\I{|\e W_1|>\tfrac{\e^\gamma}{2}, T_1< \mu_1|\ln\e|}\\
&=\I{\sup_{s\in [0,T_1]}|z^1_s(y)|> \tfrac{\e^\gamma}{2}}
+\I{\e W_1 \in I^+_{\e^\gamma}}
+\I{|\e W_1|>\tfrac{\e^\gamma}{2}, T_1< \mu_1|\ln\e|}.
\ea

\noindent
\textbf{Step 2.} Consider  $\I{x^1_s(y)\in I, s\in [0,T_1], x^1_{T_1}(y)+\e W_1\notin I}$.
For $y\in I$, we may estimate
\ba\label{eq:a2}
&\I{x^1_s(y)\in I, s\in [0,T_1], x^1_{T_1}(y)+\e W_1\notin I}\\
&\leq \I{\sup_{s\in [0,T_1]}|z^1_s(y)|> \tfrac{\e^\gamma}{2}}+
\I{\sup_{s\in [0,T_1]}|z^1_s(y)|\leq \tfrac{\e^\gamma}{2},
x^1_{T_1}(y)+\e W_1\notin I}\\
&\leq\I{\sup_{s\in [0,T_1]}|z^1_s(y)|> \tfrac{\e^\gamma}{2}}\\
&+\I{\sup_{s\in [0,T_1]}|z^1_s(y)|\leq \tfrac{\e^\gamma}{2},
|\e W_1|\leq\tfrac{\e^\gamma}{2}, T_1\geq \mu_2\e^\gamma, x^1_{T_1}(y)+\e W_1\notin I} \,\,(=0)\\
&+\I{\sup_{s\in [0,T_1]}|z^1_s(y)|\leq \tfrac{\e^\gamma}{2},
|\e W_1|\leq\tfrac{\e^\gamma}{2}, T_1< \mu_2\e^\gamma, x^1_{T_1}(y)+\e W_1\notin I}\\
&+\I{\sup_{s\in [0,T_1]}|z^1_s(y)|\leq \tfrac{\e^\gamma}{2},
|\e W_1|>\tfrac{\e^\gamma}{2}, T_1\geq \mu_1|\ln\e|, x^1_{T_1}(y)+\e W_1 \notin I}\\
&+\I{\sup_{s\in [0,T_1]}|z^1_s(y)|\leq \tfrac{\e^\gamma}{2},
|\e W_1|>\tfrac{\e^\gamma}{2},T_1< \mu_1|\ln\e|,x^1_{T_1}(y)+\e W_1\notin I }\\
&\leq\I{\sup_{s\in [0,T_1]}|z^1_s(y)|> \tfrac{\e^\gamma}{2}}
+\I{T_1< \mu_2\e^\gamma}
+\I{ \e W_1 \notin I^-_{\e^\gamma}}
+\I{|\e W_1|>\tfrac{\e^\gamma}{2},T_1< \mu_1|\ln\e|}.
\ea

\noindent
\textbf{Step 3.} Consider  $\I{x^1_s(y)\notin I\text{ for some } s\in[0,T_1]}$.
For $y\in I^-_{\e^\gamma}$, we may estimate
\ba
\label{eq:a3}
&\I{x^1_s(y)\notin I\text{ for some } s\in[0,T_1]}\leq
\I{\sup_{s\in [0,T_1]}|z^1_s(y)|> \tfrac{\e^\gamma}{2}}+\\
&\I{\sup_{s\in [0,T_1]}|z^1_s(y)|\leq \tfrac{\e^\gamma}{2},
x^1_s(y)\notin I\text{ for some } s\in[0,T_1]}(=0)
\ea

\noindent
\textbf{Step 4.} Consider  $\I{x^1_s(y)\in I, s\in(0,T_1], x^1_{T_1}(y)+\e W_1\in I,
x^2_s(x^1_{T_1}(y)+\e W_1)\notin I\text{ for some } s\in[0,T_2]}$.
For $y\in I$, we may estimate
\ba
\label{eq:a4}
&\I{x^1_s(y)\in I, s\in(0,T_1], x^1_{T_1}(y)+\e W_1\in I,
x^2_s(x^1_{T_1}(y)+\e W_1)\notin I\text{ for some } s\in[0,T_2]}\\
&=\I{x^1_s(y)\in I, s\in(0,T_1], x^1_{T_1}(y)+\e W_1\in I^-_{\e^\gamma},
x^2_s(x^1_{T_1}(y)+\e W_1)\notin I\text{ for some } s\in[0,T_2]}\\
&+\I{x^1_s(y)\in I, s\in(0,T_1], x^1_{T_1}(y)+\e W_1\in I\backslash I^-_{\e^\gamma},
x^2_s(x^1_{T_1}(y)+\e W_1)\notin I\text{ for some } s\in[0,T_2]}\\
&\leq \I{x^1_s(y)\in I, s\in(0,T_1], x^1_{T_1}(y)+\e W_1\in I^-_{\e^\gamma}}\cdot
\sup_{y\in I^-_{\e^\gamma}}\I{x^2_s(y)\notin I\text{ for some } s\in[0,T_2]}\\
&+\I{x^1_s(y)\in I, s\in(0,T_1], x^1_{T_1}(y)+\e W_1\in I\backslash I^-_{\e^\gamma}}\\
&\le \sup_{y\in I^-_{\e^\gamma}}\I{x^2_s(y)\notin I\text{ for some } s\in[0,T_2]}
+\I{x^1_s(y)\in I, s\in(0,T_1], x^1_{T_1}(y)+\e W_1\in I\backslash I^-_{\e^\gamma}}.
\ea
The first term in the resulting expression in the Step 4 is identical to the expression handled in Step 3,
while the second term requires an inessential modification of the estimation in Step 1.

Now we apply \eqref{eq:a1}, \eqref{eq:a2}, \eqref{eq:a3} and \eqref{eq:a4} to
estimate the expectations
treated in Steps 1, 2, 3 and 4 above.
In what follows, $c_\cdot$ and $C_\cdot$ denote appropriate positive constants.
Fix also $0<\delta<\min\{\gamma/2,\alpha(1-\gamma),\alpha/2\}=\min\{\gamma/2,\alpha/2\}$.

\noindent\textbf{Step 1.} Estimate
$\E\left[  \sup_{y\in I}\I{x^1_s(y)\in I, s\in(0,T_1], x^1_{T_1}+\e W_1\in I}\right].$ We get
\ba
\label{eq:esta1}
&\E\left[  \sup_{y\in I}\I{x^1_s(y)\in I, s\in(0,T_1], x^1_{T_1}+\e W_1\in I}\right]\\
&\leq \P(\sup_{[0,T_1]}|z^1_s(y)|> \tfrac{\e^\gamma}{2})+
\P(\e W_1 \in I^+_{\e^\gamma})+\P(|\e W_1|>\tfrac{\e^\gamma}{2})\P(T_1< \mu_1|\ln\e|)\\
&\leq c_1\e^{(\alpha+\gamma)/2}
+1-\frac{\e^{\alpha/2}}{2}\left[\frac{1}{(a+\e^\gamma)^\alpha}+ \frac{1}{(b+\e^\gamma)^\alpha}\right]
+  c_2 \e^{\alpha(3/2-\gamma)}|\ln\e|\\
&\leq 1-\e^{\alpha/2}\frac{\theta}{2}(1-C_1\e^\delta).
\ea

\noindent\textbf{Step 2.} Estimate
$\E \left[ \sup_{y\in I}\I{x^1_s(y)\in I, s\in(0,T_1],
x^1_{T_1}+\e W_1\notin I}\right].$ In fact,
\ba
\label{eq:esta2}
&\E \left[ \sup_{y\in I}\I{x^1_s(y)\in I, s\in(0,T_1], x^1_{T_1}+\e W_1\notin I}\right]\\
&\leq \P(\sup_{[0,T_1]}|z^1_s(y)|> \tfrac{\e^\gamma}{2})
+\P(T_1< \mu_2\e^\gamma)\\
&+\P(\e W_1 \notin I^-_{\e^\gamma})
+\P(|\e W_1|>\tfrac{\e^\gamma}{2})\P(T_1< \mu_1|\ln\e|) \\
&\leq c_1\e^{(\alpha+\gamma)/2}+c_5\e^{\alpha/2+\gamma}+
\frac{\e^{\alpha/2}}{2}\left[\frac{1}{(a-\e^\gamma)^\alpha}+ \frac{1}{(b-\e^\gamma)^\alpha}\right]
+  c_2 \e^{\alpha(3/2-\gamma)}|\ln\e|\\
&\leq\e^{\alpha/2}\frac{\theta}{2}(1+C_2\e^\delta)
\ea

\noindent\textbf{Step 3.} Estimate
$\E
\left[ \sup_{y\in I^-_{\e^\gamma}}\I{x^1_s(y)\notin I\text{ for some } s\in[0,T_1]}\right].$ We have
\ba
\label{eq:esta3}
\E \left[ \sup_{y\in I-_{\e^\gamma}}\I{x^1_s(y)\notin I\text{ for some } s\in[0,T_1]}\right]&
\leq \P(\sup_{[0,T_1]}|z^1_s(y)|> \tfrac{\e^\gamma}{2})\leq
c_1\e^{(\alpha+\gamma)/2}\leq \e^{\alpha/2}\cdot C_3\e^\delta
\ea

\noindent\textbf{Step 4.} Estimate\\
$\E\left[ \sup_{y\in I}\I{x^1_s(y)\in I, s\in(0,T_1], x^1_{T_1}(y)+\e W_1\in I,
x^2_s(x^1_{T_1}(y)+\e W_1)\notin I\text{ for some } s\in[0,T_2]}\right] $. We finally obtain
\ba
\label{eq:esta4}
&\E\left[ \sup_{y\in I}\I{x^1_s(y)\in I, s\in(0,T_1], x^1_{T_1}(y)+\e W_1\in I,
x^2_s(x^1_{T_1}(y)+\e W_1)\notin I\text{ for some } s\in[0,T_2]}\right]\\
&\leq
\P(\sup_{[0,T_1]}|z^1_s(y)|> \tfrac{\e^\gamma}{2})\\
&+ \P(\sup_{[0,T_1]}|z^1_s(y)|> \tfrac{\e^\gamma}{2})+
\P(\e W_1 \in [-b-\e^\gamma, -b+2\e^\gamma])+\P(\e W_1 \in [a-2\e^\gamma, a+\e^\gamma])\\
&+\P(|\e W_1|>\tfrac{\e^\gamma}{2})\P(T_1< \mu_1|\ln\e|)
\leq \e^{\alpha/2}\cdot C_4\e^\delta.
\ea
%%%
%%%
%%%%%%
%%%
Then for $x\in I^-_{\e^\gamma}$, $0<\e\leq \e_0$, and some positive $C_5$,
\ba
\Px(\sigma(\e)>u)&\leq
\P(\tau_1>u)
\left[ \e^{\alpha/2}\frac{\theta}{2}(1+C_2\e^\delta)+ \e^{\alpha/2}\cdot C_3\e^\delta\right] \\
&+\sum_{k=2}^\infty \P(\tau_k>u)
\left( 1-\e^{\alpha/2}\frac{\theta}{2}(1-C_1\e^\delta)\right)^{k-1}
\e^{\alpha/2}\frac{\theta}{2}
\left[ \frac{\frac{2}{\theta}C_4\e^\delta}{1-\e^{\alpha/2}\frac{\theta}{2}(1-C_1\e^\delta)}
+1+C_2\e^\delta \right] \\
&\leq \sum_{k=1}^\infty
\int_u^\infty \beta_\e e^{-\beta_\e t}\frac{(\beta_\e t)^{k-1}}{(k-1)!}\, dt
\left( 1-\e^{\alpha/2}\frac{\theta}{2}(1-C_5\e^\delta)\right)^{k-1}
\e^{\alpha/2}\frac{\theta}{2}(1+C_5\e^\delta)\\
&\leq \e^{\alpha/2}\frac{\theta}{2}(1+C_5\e^\delta)
\int_u^\infty \beta_\e e^{-\beta_\e t \e^{\alpha/2}\frac{\theta}{2}(1-C_5\e^\delta)}\,dt\\
&\leq \frac{1+C_5\e^\delta}{1-C_5\e^\delta}
\exp\left\lbrace -u \e^\alpha\frac{\theta}{\alpha}(1-C_5\e^\delta)\right\rbrace \leq
\exp\left\lbrace- u \e^\alpha\frac{\theta}{\alpha}(1-C_5\e^\delta)\right\rbrace(1+C_6\e^\delta)\\
&\leq \exp\left\lbrace -u \e^\alpha\frac{\theta}{\alpha}(1-C\e^\delta)\right\rbrace(1+C\e^\delta)
\ea
with $C=\max\{C_5, C_6\}$.

In the previous formula we have changed summation and integration. This can be done
due to the uniform convergence of the series $\sum_{k=1}^\infty
\beta_\e e^{-\beta_\e t}\frac{(\beta_\e t)^{k-1}}{(k-1)!}
\left[ 1-\e^{\alpha/2}\frac{\theta}{2}(1-C_5\e^\delta)\right]^{k-1}$ for all $t\geq 0$ and
$\e\leq \e_0$. Indeed, let $t^*_k$ be the coordinate of the maximum of the density
of  the $(\beta_\e, k)$-Gamma distribution. For $k\geq 2$, it is easy to see that
$t^*_k=\frac{k-1}{\beta_\e}$
Then, with help of Stirling's formula we obtain, that
\ba
0\leq \beta_\e e^{-\beta_\e t}\frac{(\beta_\e t)^{k-1}}{(k-1)!}\leq
\beta_\e e^{-(k-1)}\frac{(k-1)^{k-1}}{(k-1)!}
\sim \frac{1}{\sqrt{2\pi}}\frac{\beta_\e}{\sqrt{k-1}},\quad k\to \infty.
\ea
Then, for all $\e\leq \e_0$,
\ba
\sum_{k=1}^\infty
\beta_\e e^{-\beta_\e t}\frac{(\beta_\e t)^{k-1}}{(k-1)!}
\left[ 1-\e^{\alpha/2}\frac{\theta}{2}(1-C_5\e^\delta)\right]^{k-1}\leq
c_1\frac{\beta_\e}{\e^{\alpha/2}\frac{\theta}{2}(1-C_5\e^\delta)}\leq c,
\ea
where the constant $c$  does not depend on $t$ and $\e$. Thus uniform convergence
follows from dominated convergence.

\end{proof}

\subsection{Lower estimate \label{S:Si:below}}

In this subsection we estimate $\Px(\sigma(\e)>u)$  from below as $\e\to 0$, $u>0$. This leads to
the following Proposition with a rather technical proof again.

\begin{prop}
\label{p:below}
There exist constants $\e_0>0$ and $C>0$ such that for all $0<\e\leq \e_0$,
$0<\delta<\min\{\alpha/2,\gamma/2\}$, $x\in [-b+\e^\gamma, a-\e^\gamma]$ and $u\geq 0$ the estimate
\begin{equation}
\Px(\sigma(\e)>u)\geq
\exp\left\lbrace
-u \frac{\e^\alpha}{\alpha}\left[\frac{1}{a^\alpha}+\frac{1}{b^\alpha}\right] (1+C\e^\delta)
\right\rbrace
(1-C\e^\delta)
\end{equation}
is valid.

\end{prop}
\begin{proof}
We use the following inequality:
\begin{equation}
\begin{aligned}
\Px(\sigma(\e)>u)&=\sum_{k=1}^\infty \P(\tau_k>u)\Px(\sigma(\e)=\tau_k)+
\Px(\sigma(\e)>u|\sigma\in (\tau_{k-1},\tau_k))\Px(\sigma\in (\tau_{k-1},\tau_k))\\
&\geq  \sum_{k=1}^\infty \P(\tau_k>u) \Px(\sigma(\e)=\tau_k).
\end{aligned}
\end{equation}

With arguments analogous to \eqref{eq:stauk} we obtain the factorisation
\begin{equation}
\begin{aligned}
&\Px(\sigma=\tau_k)=\E_x\I{X_s^\e\in I, s\in[0,\tau_k), X_{\tau_k}^\e\notin I}\\
&= \E_x \prod_{j=1}^{k-1}
\I{x^j_s(x^{j-1}_{T_{j-1}}+\e W_{j-1})\in I, s\in[0,T_j],
x^j_{T_j}(x^{j-1}_{T_{j-1}}+\e W_{j-1})+\e W_j\in I} \\
&\qquad\qquad\qquad\times
\I{x^k_s(x^{k-1}_{T_{k-1}}+\e W_{k-1})\in I, s\in[0,T_k],
x^k_{T_k}(x^{k-1}_{T_{k-1}}+\e W_{k-1})+\e W_k\notin I}\\
&\geq \E_x\prod_{j=1}^{k-1}
\I{x^j_s(x^{j-1}_{T_{j-1}}+\e W_{j-1})\in I, s\in[0,T_j],
x^j_{T_j}(x^{j-1}_{T_{j-1}}+\e W_{j-1})+\e W_j\in I^-_{\e^\gamma}} \\
&\qquad\qquad\qquad\times
\I{x^k_s(x^{k-1}_{T_{k-1}}+\e W_{k-1})\in I, s\in[0,T_k],
x^k_{T_k}(x^{k-1}_{T_{k-1}}+\e W_{k-1})+\e W_k\notin I}\\
&\geq \E\prod_{j=1}^{k-1}
\inf_{y\in I^-_{\e^\gamma}}\I{x^j_s(y)\in I, s\in(0,T_j],
x^j_{T_j}(y)+\e W_j\in I^-_{\e^\gamma}} \\
&\qquad\qquad\qquad\times
\inf_{y\in I^-_{\e^\gamma}}\I{x^k_s(y)\in I, s\in(0,T_k], x^k_{T_k}(y)+\e W_k\notin I}\\
&=\prod_{j=1}^{k-1}
\E\left[ \inf_{y\in I^-_{\e^\gamma}}\I{x^j_s(y)\in I, s\in(0,T_j],
x^j_{T_j}(y)+\e W_j\in I^-_{\e^\gamma}}\right] \\
&\qquad\qquad\qquad\times
\E \left[ \inf_{y\in I^-_{\e^\gamma}}\I{x^k_s(y)\in I, s\in(0,T_k],
x^k_{T_k}(y)+\e W_k\notin I}\right] \\
=&\left( \E\left[\inf_{y\in I^-_{\e^\gamma}}\I{x^1_s(y)\in I, s\in(0,T_1],
x^1_{T_1}(y)+\e W_1\in I^-_{\e^\gamma}}\right]\right)^{k-1}  \\
&\qquad\qquad\qquad\times
\E \left[\inf_{y\in I^-_{\e^\gamma}}\I{x^1_s(y)\in I, s\in(0,T_1],
x^1_{T_1}(y)+\e W_1\notin I}\right].
\end{aligned}
\end{equation}
For $y\in I^-_{\e^{\gamma}}$, we next specify separately in two steps the further estimation for the two different events appearing in the formulae for
$\Px(\sigma(\e)=\tau_k)$.

\noindent
\textbf{Step 1.}
Consider the event $\I{x^1_s(y)\in I, s\in (0,T_1], x^1_{T_1}(y)+\e W_1\in I^-_{\e^\gamma}}.$ We may estimate
\ba
\label{eq:b1}
&\I{x^1_s(y)\in I, s\in (0,T_1], x^1_{T_1}(y)+\e W_1\in I^-_{\e^\gamma}}\\
&\geq
\I{\sup_{s\in [0,T_1]}|z^1_s(y)|\leq \tfrac{\e^\gamma}{2},
 x^1_s(y)\in I, s\in (0,T_1], x^1_{T_1}(y)+\e W_1\in I^-_{\e^\gamma}}\\
&\geq
\I{\sup_{s\in [0,T_1]}|z^1_s(y)|\leq \tfrac{\e^\gamma}{2},
 x^1_s(y)\in I, s\in (0,T_1],|\e W_1|\leq\tfrac{\e^\gamma}{2}, T_1\geq \mu_2\e^\gamma,
x^1_{T_1}(y)+\e W_1\in I^-_{\e^\gamma}}\\
&+\I{\sup_{s\in [0,T_1]}|z^1_s(y)|\leq \tfrac{\e^\gamma}{2},
 x^1_s(y)\in I, s\in (0,T_1],|\e W_1|>\tfrac{\e^\gamma}{2}, T_1\geq \mu_1|\ln\e|,
x^1_{T_1}(y)+\e W_1\in I^-_{\e^\gamma}}\\
&\geq
\I{\sup_{s\in [0,T_1]}|z^1_s(y)|\leq \tfrac{\e^\gamma}{2},|\e W_1|\leq\tfrac{\e^\gamma}{2},
T_1\geq \mu_2\e^\gamma}\\
&+\I{\sup_{s\in [0,T_1]}|z^1_s(y)|\leq \tfrac{\e^\gamma}{2},
|\e W_1|>\tfrac{\e^\gamma}{2}, T_1\geq\mu_1|\ln\e|, \e W_1\in I^-_{2\e^\gamma}}\\
&=\I{|\e W_1|\leq\tfrac{\e^\gamma}{2}, T_1\geq \mu_2\e^\gamma}\\
&-\I{\sup_{s\in [0,T_1]}|z^1_s(y)|>\tfrac{\e^\gamma}{2},|\e W_1|\leq\tfrac{\e^\gamma}{2},
T_1\geq \mu_2\e^\gamma}\\
&+\I{|\e W_1|>\tfrac{\e^\gamma}{2}, T_1\geq \mu_1|\ln\e|, \e W_1\in I^-_{2\e^\gamma}}\\
&-\I{\sup_{s\in [0,T_1]}|z^1_s(y)|> \tfrac{\e^\gamma}{2},
|\e W_1|>\tfrac{\e^\gamma}{2}, T_1\geq \mu_1|\ln\e|, \e W_1\in I^-_{2\e^\gamma}}\\
&=\I{|\e W_1|\leq\tfrac{\e^\gamma}{2}}-
\I{|\e W_1|\leq\tfrac{\e^\gamma}{2}, T_1< \mu_2\e^\gamma}
-\I{\sup_{s\in [0,T_1]}|z^1_s(y)|>\tfrac{\e^\gamma}{2},|\e W_1|\leq\tfrac{\e^\gamma}{2},
T_1\geq \mu_2\e^\gamma}\\
&+\I{|\e W_1|>\tfrac{\e^\gamma}{2}, \e W_1\in I^-_{2\e^\gamma}}-
\I{|\e W_1|>\tfrac{\e^\gamma}{2}, T_1< \mu_1|\ln\e|, \e W_1\in I^-_{2\e^\gamma}}
\\
&-\I{\sup_{s\in [0,T_1]}|z^1_s(y)|> \tfrac{\e^\gamma}{2},
|\e W_1|>\tfrac{\e^\gamma}{2}, T_1\geq \mu_1|\ln\e|, \e W_1\in I^-_{2\e^\gamma}}\\
&\geq \I{\e W_1\in I^-_{2\e^\gamma}}-
\I{T_1< \mu_2\e^\gamma}-2\I{\sup_{s\in [0,T_1]}|z^1_s(y)|>\tfrac{\e^\gamma}{2}}
-\I{|\e W_1|>\tfrac{\e^\gamma}{2}, T_1< \mu_1|\ln\e|}.
\ea

\noindent
\textbf{Step 2.}
The event $\I{x^1_s(y)\in I, s\in(0,T_1], x^1_{T_1}(y)+\e W_1\notin I}$ may be estimated as follows:
\ba
\label{eq:b2}
&\I{x^1_s(y)\in I, s\in(0,T_1], x^1_{T_1}(y)+\e W_1\notin I}\\
&\geq  \I{\sup_{s\in [0,T_1]}|z^1_s(y)|\leq \tfrac{\e^\gamma}{2},
x^1_s(y)\in I, s\in(0,T_1], x^1_{T_1}(y)+\e W_1\notin I}\\
&\geq \I{\sup_{s\in [0,T_1]}|z^1_s(y)|\leq \tfrac{\e^\gamma}{2},
x^1_s(y)\in I, s\in(0,T_1], T_1\geq \mu_1|\ln\e|, x^1_{T_1}(y)+\e W_1\notin I}\\
&\geq \I{\sup_{s\in [0,T_1]}|z^1_s(y)|\leq \tfrac{\e^\gamma}{2},
T_1\geq \mu_1|\ln\e|, \e W_1\notin I^+_{\e^\gamma}}\\
&=\I{T_1\geq \mu_1|\ln\e|, \e W_1\notin I^+_{\e^\gamma}}
-\I{\sup_{s\in [0,T_1]}|z^1_s(y)|> \tfrac{\e^\gamma}{2},
T_1\geq \mu_1|\ln\e|, \e W_1\notin I^+_{\e^\gamma}}\\
&\geq \I{\e W_1\notin I^+_{\e^\gamma}}-\I{T_1< \mu_1|\ln\e|, \e W_1\notin I^+_{\e^\gamma}}
-\I{\sup_{s\in [0,T_1]}|z^1_s(y)|> \tfrac{\e^\gamma}{2}, \e W_1\notin I^+_{\e^\gamma}}.
\ea
Now we apply \eqref{eq:b1} and \eqref{eq:b2}  to estimate the expectations
appearing in the formula for $\Px(\sigma=\tau_k)$.
In what follows, $c_i$ and $C_i$ denote appropriate positive constants.
Fix also $0<\delta<\min\{\gamma/2,\alpha(1-\gamma),\alpha/2\}=\min\{\gamma/2,\alpha/2\}$.

\noindent\textbf{Step 1.} Here we estimate
$\E\left[ \inf_{y\in I^-_{\e^\gamma}}\I{x^1_s(y)\in I, s\in(0,T_1],
x^1_{T_1}+\e W_1\in I^-_{\e^\gamma}}\right].$ In fact, we obtain from employing results from sections
\ref{s:dev} and \ref{s:heur}

\ba
&\E\left[ \inf_{y\in I^-_{\e^\gamma}}\I{x^1_s(y)\in I, s\in(0,T_1],
x^1_{T_1}+\e W_1\in I^-_{\e^\gamma}}\right]\\
&\geq  \P(\e W_1\in I^-_{2\e^\gamma})-
\P(T_1< \mu_2\e^\gamma)-2\P(\sup_{[0,T_1]}|z^1_s(y)|>\tfrac{\e^\gamma}{2})
-\P(|\e W_1|>\tfrac{\e^\gamma}{2}, T_1< \mu_1|\ln\e|)\\
&\geq 1-\frac{\e^{\alpha/2}}{2}\left[\frac{1}{(a-2\e^\gamma)^\alpha}+
\frac{1}{(b-2\e^\gamma)^\alpha}\right]- c_1 \e^{\alpha/2+\gamma}-2c_2\e^{(\alpha+\gamma)/2}
-c_3\e^{\alpha(3/2-\gamma)}|\ln\e|\\
&\geq 1-\e^{\alpha/2}\frac{\theta}{2}(1+C_1\e^\delta).
\ea

\noindent\textbf{Step 2.} We next estimate
$\E\left[\inf_{y\in I^-_{\e^\gamma}}
\I{x^1_s(y)\in I, s\in(0,T_1], x^k_{T_1}+\e W_1\notin I}\right],$ for which we obtain similarly

\ba
&\E\left[\inf_{y\in I^-_{\e^\gamma}}
\I{x^1_s(y)\in I, s\in(0,T_1], x^k_{T_1}+\e W_1\notin I}\right]\geq \\
&\P(\e W_1\notin I^+_{\e^\gamma})
\left(1 -\P(T_1< \mu_1|\ln\e|)
-\P(\sup_{[0,T_1]}|z^1_s(y)|> \tfrac{\e^\gamma}{2})\right)\\
&\geq  \frac{\e^{\alpha/2}}{2}\left[\frac{1}{(a+\e^\gamma)^\alpha}+
 \frac{1}{(b+\e^\gamma)^\alpha}\right]\left( 1-
c_3\e^{\alpha/2}|\ln\e|-c_2\e^{(\alpha+\gamma)/2}\right) \\
&\geq  \e^{\alpha/2}\frac{\theta}{2}\left( 1-C_2\e^{\delta}\right).
\ea
Consequently, with $C_3=\max\{C_1,C_2\}$, for $x\in I^-_{\e^\gamma}$, $0<\e\leq \e_0$
\ba
\Px(\sigma(\e)>u)&\geq \sum_{k=1}^\infty \int_u^\infty
\beta e^{-\beta t}\frac{(\beta t)^{k-1}}{(k-1)!}\, dt
\left[ 1-  \e^{\alpha/2}\frac{\theta}{2}  (1+C_3\e^\delta)\right]^{k-1}
\e^{\alpha/2}\frac{\theta}{2}(1-C_3\e^\delta) \\
&=\e^{\alpha/2}\frac{\theta}{2}(1-C_4\e^\delta)
\int_u^\infty \beta e^{-\beta t \e^{\alpha/2}\frac{\theta}{2}(1+C_3\e^\delta)}\,dt\\
&\geq \frac{1-C_4\e^\delta}{1+C_3\e^\delta}
\exp\left\lbrace-u \e^\alpha\frac{\theta}{\alpha}(1+C_3\e^\delta)\right\rbrace\geq
\exp\left\lbrace -u \e^\alpha\frac{\theta}{\alpha}(1+C_3\e^\delta)\right\rbrace (1-C_5\e^\delta)\\
&\geq \exp\left\lbrace -u \e^\alpha\frac{\theta}{\alpha}(1+C\e^\delta)\right\rbrace (1-C\e^\delta)
\ea
with $C=\max\{C_3,C_5\}$.
See the end of the proof of Proposition \ref{p:above} for the justification of switching the order
of summation and integration in the above argument.

\end{proof}

\begin{Aproof}{Theorem}{\ref{th:main}}
The statement of Theorem \ref{th:main} follows directly from Propositions \ref{p:above}
 and \ref{p:below}.

\end{Aproof}

\begin{rem}
\label{rm:1} {\rm

The threshold $1/\!\sqrt{\e}$ for separating the two parts of $L$ is not the only possible choice.
Indeed, let us consider thresholds of the type $1/\e^\rho$ for $\rho>0$.
Then $\beta_\e=\int_{\mathbb{R}}\I{|y|>\frac{1}{\e^\rho}}\nu(dy)=\frac{2}{\alpha}\e^{\alpha\rho}$.
$\rho$ must satisfy some further conditions. Firstly, we demand that $\rho<1$ so that
we can easily calculate the probability
\be
P(\e W_1\notin [-b,a])=
\frac{\alpha}{2\e^{\alpha\rho}}\int_{\mathbb{R}\backslash [-\frac{b}{\e},\frac{a}{\e}]}
\I{|y|>\frac{1}{\e^\rho}}\frac{dy}{|y|^{1+\alpha}}=\mathcal{O}(\e^{\alpha(1-\rho)}).
\ee
This is the \emph{characteristic} probability of our analysis. The probabilities of other relevant
events should have smaller order in the small noise limit. This for example applies to the event that
$\xi^\e$ leaves the $\e^\gamma$-tube around the deterministic trajectory (see Proposition
\ref{p:1}), and  the event that $T_1<\mu_2\e^\gamma)$
(e.g.\ see \eqref{eq:esta1}, \eqref{eq:esta2} and \eqref{eq:esta3}).
This leads to the following inequalities on $\rho$ and $\gamma$:
$0<\rho<1$, $\gamma>0$, $\alpha(1-\rho)<2-2\rho-2\gamma$ and $\alpha(1-\rho)<\alpha\rho+\gamma$.
Applying some algebra we rewrite these inequalities in the form
\be
\begin{cases}
\gamma<\frac{2-\alpha}{2}(1-\rho),\qquad 0<\rho<1, \\
\gamma>\alpha(1-2\rho), \qquad \qquad\gamma>0.
\end{cases}
\ee
The solution set $(\rho,\gamma)$ is non-empty for all $\alpha\in (0,2)$ and depends on $\alpha$.
However $\rho=\frac{1}{2}$ is the minimal value independent of $\alpha$ for which there exist
$\gamma$ solving the inequalities. In this case any $\gamma$ from the interval
$(0,\frac{2-\alpha}{4})$ is a solution. For our purposes we have taken
$\gamma=\frac{2-\alpha}{5}$.

}
\end{rem}

\section{Return from $-\infty$ and deviations
from the deterministic trajectory: exit from unbounded interval \label{s:devU}}

With the aim of proving an analogue of Proposition~\ref{p:1}, we study in this section
the exit problem of the solution of our stochastic differential equation from the unbounded
interval $J=(-\infty, a]$.
In this case we shall use the condition that $U'$ increases faster than a linear function at
$-\infty$ which guarantees a return from infinity in finite
time for the unperturbed deterministic motion.
For simplicity we assume that for large $|x|$, $U$ is a power function, i.e.\
$|U(x)|= c_1|x|^{2+c_2}$, $c_1,c_2>0$ for $|x|\ge N$.
This condition can be weakened (see Remark~\ref{rm:2} at the end of this section).
We stick to it to avoid
technicalities irrelevant for the main aim of the paper.

Fix two more positive numbers $r$ and $R$ such that $N<r<R$, and such that with
$T_R=\int_{-\infty}^{-R+10}\frac{dy}{|U'(y)|}$ we have $-r=Y_{T_R}(-R)$.
Consider the equations \eqref{eq:eqs} on the unbounded interval $J$. As in section \ref{s:dev}
we estimate on the basis of the representations $x^\e(x)=Y_t(x)+\e Z^\e_t(x)+R^\e_t(x)$, with
$Y_t(x)=x-\int_0^t U'(Y_s(x))\, ds$ and
$Z^\e_t(x)=\xi^\e_t-\int_0^t \xi^\e_{s-}U''(Y_s(x))\frac{U'(Y_t(x))}{U'(Y_s(x))}\,ds.$
We first prove an estimate enabling us to transfer Lemma \ref{l:Z} to unbounded intervals.

\begin{lemma}
\label{l:Zinf}
The inequality
$\sup_{t\in [0,T_R]}|Z^\e_t(x)|\leq 2\sup_{t\in [0,T_R]}|\xi^\e_t|$
holds a.s.\ for $x\leq -R$.
\end{lemma}
\begin{proof}
For all $x\leq -R$ by definition of $R$ and $r$ we have $Y_t(x)\leq -r$. Moreover,
by assumption $U''(Y_t(x))>0$ for $t\in[0,T_R]$, whence
\be
Z^\e_t(x)\leq\sup_{t\in [0,T_R]}|\xi^\e_t|
\left(1+\sup_{x\leq -R}\sup_{t\in [0,T_R]}
\int_0^t U''(Y_s(x))\frac{U'(Y_t(x))}{U'(Y_s(x))}\,ds\right).
\ee
We show that the integral in the latter parenthesis is uniformly bounded in $x$.
Denote $Y_t(x)=v$. Then $dv=\dot Y_t(x)\, dt=-U'(Y_t(x))\, dt=-U'(v)\,dt$. Therefore
\begin{eqnarray*}
\int_0^t U''(Y_s(x))\frac{U'(Y_t(x))}{U'(Y_s(x))}\,ds &=&
-U'(Y_t(x))\int_x^{Y_t(x)} \frac{U''(v)}{U'(v)^2}\,dv\\ &=&
U'(Y_t(x))\left[ \frac{1}{U'(Y_t(x))}-\frac{1}{U'(x)}\right] \leq 1.
\end{eqnarray*}

\end{proof}

Here is the analogue of Lemma \ref{l:Z}.

\begin{lemma}
\label{l:ZinfR}
There is a constant $C'_Z>0$ such that the inequality
\be
\sup_{t\in [0,T]}|Z^\e_t(x)|\leq C'_Z\sup_{t\in [0,T]}|\xi^\e_t|
\ee
holds a.s.\ for $x\leq a$, $T > 0, \e>0$.
\end{lemma}

\begin{proof}
The proof obviously has to combine the previous Lemma with Lemma
\ref{l:Z} on bounded intervals. The inequality holds with $C'_Z=1+C_Z$.
\end{proof}

To estimate deviations of the random paths from the paths of the deterministic equation, we
restrict from the start to sets of bounded scaled noise.
\begin{lemma}
\label{l:xYinf}
On the event $\{ \sup_{t\in [0,T_R]}|\e\xi^\e_t|  < 1\}$ the
following inequality holds a.s.\
\be
\sup_{t\in [0,T_R]} |x^\e_t-Y_t(x)| < 10
\ee
uniformly for $x\leq -R$.
\end{lemma}
\begin{proof}
It follows from Lemma~\ref{l:Zinf}, that
$\sup_{t\in [0,T_R]}|\e Z_t^\e(x)| < 2 $ a.s.\ on the event
$\{ \sup_{t\in [0,T_R]}|\e\xi^\e_t| < 1 \}$.
We show that the rest term $|R^\e|$ is bounded by $8$. Indeed, the rest term satisfies
the integral equation
\be
\label{eq:Rint}
R^\e_t(x)=\int_0^t \left[ -U'(Y_s(x)+\e Z_{s-}^\e(x)+R_s^\e(x))+ U'(Y_s(x))+
U''(Y_s(x))\e Z_{s-}^\e(x) \right]\, ds.
\ee
$R^\e(x)$ is absolutely continuous a.s.\ and   $R_0^\e(x)=0$.
Assume, that there exists a smallest $\tau\in [0,T_R]$ such that $R^\e_\tau(x)=8$. Then
the left Dini derivative of $R^\e(x)$ at this point is non-negative, i.e.\
\be
-U'(Y_\tau(x)+\e Z_{\tau-}^\e(x)+8)+ U'(Y_\tau(x))+
U''(Y_\tau(x))\e Z_{\tau-}^\e(x)\geq 0.
\ee
On the other hand, our conditions on $U$ guarantee
\ba
&-U'(Y_\tau(x)+\e Z_{\tau-}^\e(x)+8)+ U'(Y_\tau(x))+
U''(Y_\tau(x))\e Z_{\tau-}^\e(x)\\
&\qquad\qquad
<-U'(Y_\tau(x)+6)+ U'(Y_\tau(x))+ 2 U''(Y_\tau(x))<0,
\ea
and a contradiction is reached. The estimate $R^\e_t(x)>-8$ is obtained analogously.
\end{proof}

The following Lemma is an analogue of Lemma~\ref{l:Rrough} and gives a rough estimate for the
remainder term $R^\e$.

\begin{lemma}
\label{l:Rroughinf}
There exists $C>0$ such that for $x\in[-R,a]$, $T>0$
\be
\sup_{t\in[0,T]}|R^\e_t|<C
\ee
a.s.\ on the event $\{\sup_{t\in [0,T]}|\e \xi_t^\e| < 1\}.$
\end{lemma}

Lemma \ref{l:Rroughinf} again has localising consequences: It states precisely that the solution process
$x^\e$, with initial state confined to $[-R,a],$ stays bounded by a deterministic constant on sets
of the form $\{ \sup_{t\in [0,T]}|\e\xi^\e_t(\omega)| <1 \}.$ Therefore, in the small noise limit,
only local properties of $U$ are relevant to our analysis.

Let us paraphrase the most important aspect of what we found in the previous
Lemmas due to finite return from $-\infty$.
With probability close to one, the random trajectory starting at $x\leq -R$ reaches the
finite interval $[-R,-r]$ in a finite non-random time $T_R$ which does not depend on $\e$.
Our investigation therefore reduces to the study of the dynamics
of paths starting in the finite interval $[-R,a]$.
Since the deterministic trajectories starting in $[-R,a]$ do not leave this interval,
the statement of Lemma~\ref{l:Z} does not change.
Due to Lemma~\ref{l:Rroughinf}
the estimate of the rest term $R^\e$ given in
Lemma \ref{l:R} also holds unchanged.

Thus, in the case of the unbounded interval we have all necessary tools to estimate
the exit probabilities.

\begin{rem}
\label{rm:2}
  {\rm
The conditions on the behaviour of the potential $U$ at $-\infty$ can be weakened. Indeed,
a slight extension of the proof of Lemma~\ref{l:Zinf} allows to drop the convexity condition
$U''>0$. Furthermore, to show that $R^\e$ is bounded from above by some $p_1>0$,
we need to guarantee that the integrand in
\eqref{eq:Rint} is negative for $R_\tau^\e=p_1$ under the condition
$|\e Z_{\tau -}^\e|\leq p_2$, $p_2>0$. This leads to the inequality
$-\inf_{\theta\in[p_2-p_1,p_2+p_1]}U'(y+\theta)+U'(y)+p_1|U''(y)|<0$, which has to hold
for $y\leq -N$. For instance, for the power function $U(x)=c_1|x|^{2+c_2}$ considered above this
inequality is equivalent to $-(p_2-2p_1)|y|^{-1}+\mathcal{O}(|y|^{-2})<0$, $y\to-\infty$,
and thus holds for any $p_2>2p_1>0$.

}
\end{rem}

\section{The law of $\tau(\e)$ \label{s:law:tau}}

In this section we estimate $\Px(\tau(\e)>u)$ for $u\geq 0$ as $\e\to 0$.
Indeed, the extensions of Propositions \ref{p:above} for the estimation above
and \ref{p:below} for lower bounds with
parameter $b=+\infty$ take the following form.

\begin{prop}
%\label{p:above}
Let $\delta=\min\{\alpha/2,\gamma/2\}$. There exists $\e_0>0$ and $C>0$ such that for all
$0<\e\leq \e_0$, $x\in (-\infty, a-\e^\gamma]$ and $u\geq 0$
\be
\Px(\tau(\e)>u)\leq
\exp\left\lbrace
-u \frac{\e^\alpha}{\alpha a^\alpha}(1-C\e^\delta)
\right\rbrace
(1+C\e^\delta).
\ee
\end{prop}

\begin{prop}
%\label{p:above}
Let $\delta<\min\{\alpha/2,\gamma/2\}$. There exists $\e_0>0$ and $C>0$ such that for all
$0<\e\leq \e_0$, $x\in (-\infty, a-\e^\gamma]$ and $u\geq 0$
\be
\Px(\tau(\e)>u)\leq
\exp\left\lbrace
-u \frac{\e^\alpha}{\alpha a^\alpha}(1+C\e^\delta)
\right\rbrace
(1-C\e^\delta).
\ee

\end{prop}

\begin{proof}
The arguments to prove the estimates are similar to those of the bounded
case. We just need to adapt Steps 1, 2, 3 and 4
(see sections \ref{S:Si:above} and \ref{S:Si:below}) to the case of an unbounded interval.
Let us
consider for example the extension of Step 1 from section \ref{S:Si:above}.
The basic formula \eqref{eq:stauk} holds in the case of unbounded intervals with $J$ replacing $I$.
We demonstrate how to modify the reasoning just in the series of inequalities \eqref{eq:a1} and in
the estimate \eqref{eq:esta1}.
The other estimates are obtained analogously.

Firstly, we estimate $\I{x^1_s(y)\in J, s\in [0,T_1], x^1_{T_1}(y)+\e W_1\in J}$
for $-R\leq y\leq a$. Denote $A=\{ \sup_{t\in [0,T_1]}|\e\xi^\e_t|<1 \}$.
On the event $A$, the trajectory $x^1_t(y)$, $t\in[0,T_1]$, belongs to a compact interval,
so its dynamics is
indistinguishable form the one treated in the bounded case. Therefore, we have
\ba
&\I{x^1_s(y)\in J, s\in [0,T_1], x^1_{T_1}(y)+\e W_1\in J}\\
&=\I{x^1_s(y)\in J, s\in [0,T_1], x^1_{T_1}(y)+\e W_1\in J}
\left(\I{A}+\I{A^c} \right) \\
&\leq
\I{x^1_s(y)\in J, s\in [0,T_1], x^1_{T_1}(y)+\e W_1\in J,A}+\I{A^c}\leq \cdots\\
&\leq
\I{\sup_{[0,T_1]}|x^1_s(y)-Y_s(y)|> \tfrac{\e^\gamma}{2}}
+\I{\e W_1 \in J^+_{\e^\gamma}}
+\I{|\e W_1|>\tfrac{\e^\gamma}{2}, T_1< \mu_1|\ln\e|}+\I{A^c}
\ea

The case $y\leq -R$ is slightly more complicated since we have to treat the return of
$x^1(y)$ to a compact interval in a finite time.
Denote
\ba
B=\{\omega \,:\,  \sup_{[0,T_R\land T_1]}|\e\xi^\e_t|<1  \}&\cap
\{\omega \,:\,  \sup_{[T_R\land T_1, T_1]}|\e\xi^\e_t-\e\xi^\e_{T_R\land T_1}|<1\}\\
&\supseteq \{\omega \,:\,  \sup_{[0, T_1]}|\e\xi^\e_t|<1 \}=A.
\ea
Then we have
\ba
&\I{x^1_s(y)\in J, s\in [0,T_1], x^1_{T_1}(y)+\e W_1\in J}\\
&\leq \I{x^1_s(y)\in J, s\in [0,T_R\land T_1],
x^1_s(y)\in J, s\in [T_R\land T_1,T_1], x^1_{T_1}(y)+\e W_1\in J} \cap B +\I{A^c} \\
&\leq
\I{x^1_s(y)\in J, s\in [T_R\land T_1,T_1], x^1_{T_1}(y)+\e W_1\in J} \cap B+\I{A^c} \leq \cdots\\
&\leq
\I{\sup_{s\in [T_R\land T_1,T_1]}|x^1_s(y)-Y_{s-T_R\land T_1}(x^1_{T_R\land T_1}(y))|>
\tfrac{\e^\gamma}{2}}
+\I{\e W_1 \in J^+_{\e^\gamma}}\\
&+\I{|\e W_1|>\tfrac{\e^\gamma}{2}, T_R\land T_1< \mu_1|\ln\e|}+\I{A^c}\\
&\leq
\sup_{y\in [-R,-r]}\I{\sup_{s\in [0,T_1]}|x^1_s(y)-Y_s(y)|> \tfrac{\e^\gamma}{2}}
+\I{\e W_1 \in J^+_{\e^\gamma}}\\
&+\I{|\e W_1|>\tfrac{\e^\gamma}{2}, T_1< \mu_1|\ln\e|}+\I{A^c}.
\ea
These estimates may be treated in a way similar to \eqref{eq:esta1}. In fact,
\ba
\E\left[  \sup_{y\in J}\I{x^1_s(y)\in J, s\in(0,T_1], x^1_{T_1}+\e W_1\in J}\right]
&\leq 1-\e^{\alpha/2}\frac{1}{2a^\alpha}(1-C\e^\delta)
\ea
for some $C>0$.
The other steps are modified analogously.
\end{proof}

\vspace{.5cm}

\begin{Aproof}{Theorem}{\ref{th:mainU}}
Combine the estimates of the above Propositions.
\end{Aproof}

%%%%%\v{C}eby\v{s}\"ev

\bibliography{biblio-new}
\bibliographystyle{alpha}

\vspace{1cm}

\parbox{.48\linewidth}{
\noindent
Peter Imkeller\\
Institut f\"ur Mathematik\\
Humboldt Universit\"at zu Berlin\\
Rudower Chaussee 25\\
12489 Berlin Germany\\
E.mail: imkeller@mathematik.hu-berlin.de\\
http://www.mathematik.hu-berlin.de/\~{}imkeller
}\hfill
\parbox{.48\linewidth}{
\noindent
Ilya Pavlyukevich\\
Institut f\"ur Mathematik\\
Humboldt Universit\"at zu Berlin\\
Rudower Chaussee 25\\
12489 Berlin Germany\\
E.mail: pavljuke@mathematik.hu-berlin.de\\
http://www.mathematik.hu-berlin.de/\~{}pavljuke
}

\end{document}